\documentclass[letterpaper,12pt]{article}

\usepackage{amssymb,amsthm,amsmath,wasysym}
\usepackage{enumitem}
\numberwithin{equation}{section}
\usepackage{color,bbold,soul,txfonts}
\usepackage[left=0.65in,right=0.65in,top=1.00in,bottom=1.00in,footskip=0.7in]{geometry}
\usepackage[pdfpagelabels]{hyperref}
\usepackage{arydshln}
\usepackage[titletoc]{appendix}
\usepackage[nottoc,numbib]{tocbibind}

\usepackage{graphicx,subfig,float}
\newtheorem{theorem}{Theorem}[subsection]
\newtheorem{definition}[theorem]{Definition}

\newtheorem{proposition}[theorem]{Proposition}
\newtheorem{lemma}[theorem]{Lemma}
\newtheorem{corollary}[theorem]{Corollary}

\title{Equivariant Solutions to a System of Nonlinear Wave Equations with Ginzburg-Landau Type Potential}
\author{Kyle Thompson}
\date{}

\begin{document}

\maketitle

\abstract{It is known that there exist solutions with interfaces to various scalar nonlinear wave equations. In this paper, we look for solutions of a two-component system of nonlinear wave equations where one of the components has an interface and and where the second component is exponentially small except near the interface of the first component. A formal asymptotic expansion suggests that there exist solutions to this system with these characteristics whose profiles are determined by the winding number density of the second component and where the interface of the first component is a time-like surface in Minkowski space whose geometric evolution is coupled in a highly nonlinear way to the phase of the second component. We verify this heuristic when $n=2$ and for equivariant maps.}

\section{Introduction}
\subsection{Synopsis}\label{synopsis}
In this paper we consider two-component systems of hyperbolic system of PDEs qualitatively similar to
\begin{equation}\label{equations of motion of interface with a current}
 \left\{
  \begin{array}{c}
   \partial_{tt}\phi - \Delta \phi + \frac{\lambda_\phi}{\epsilon^2} (\phi^2 - 1)\phi = -\frac{\beta}{\epsilon^2}\left|\sigma\right|^2\phi \\
   \partial_{tt}\sigma - \Delta \sigma + \frac{\lambda_\sigma}{\epsilon^2} (\left|\sigma\right|^2 - 1)\sigma = -\frac{\beta}{\epsilon^2}\phi^2\sigma
  \end{array}
 \right.
\end{equation}
where $\Phi := (\phi,\sigma):\mathbb{R}^{1+n}\rightarrow \mathbb{R}\times\mathbb{C}$, $0 < \epsilon \ll 1$ is a small parameter of the model, and $(\lambda_\phi,\lambda_\sigma,\beta)$ are real, non-negative constants. We are interested in solutions to (\ref{equations of motion of interface with a current}) with the properties that
\begin{itemize}
 \item $\phi$ has an interface
 \item $\sigma$ is exponentially small except near the interface
\end{itemize}
For the first equation of (\ref{equations of motion of interface with a current}), if the right hand side vanishes (which happens if $\beta = 0$ or if $\sigma = 0$), then it is known that there exists a $\phi$ with an interface solving this equation \cite{jerrard2011semilinear}. We, however, are interested in regimes where $\phi$ and $\sigma$ are coupled (i.e. $\beta \neq 0$) and where $(\lambda_\phi,\lambda_\sigma,\beta)$ are chosen so that $(\phi,\sigma)$ have the properties described above, which in particular stipulate that $\sigma \neq 0$ near the interface of $\phi$. For these regimes, it follows from the physics literature on superconducting strings, reviewed in section \ref{physics background} below, that the $\sigma$-field can naturally be identified with a superconducting current confined to the interface of $\phi$. Hence, we call (\ref{equations of motion of interface with a current}) the \textbf{superconducting interface model}. The goal of this paper is to understand the coupling between the 
current and the interface and, in particular, how the current affects the dynamics of the interface.
\bigskip

As discussed in appendix \ref{formal asymptotics}, a formal asymptotic expansion suggests that for suitable local coordinates $(y^\tau,y^\nu) = (y_0,..,y_n)$ near a codimension one time-like surface $\Gamma$, with $y^\tau = (y_0,..,y_{n-1})$ parameterizing $\Gamma$ and with $\left\{y^\nu = 0\right\}$ corresponding to $\Gamma$, then there \textit{should} exist a solution to (\ref{equations of motion of interface with a current}) satisfying
\begin{equation}\label{formal asymp intro}
 \left\{
  \begin{array}{l}
   \phi(y^\tau,y^\nu) \approx \phi_0(\frac{y^\nu}{\epsilon};\zeta(y^\tau)) \\
   \sigma(y^\tau,y^\nu) \approx e^{\frac{i}{\epsilon}\theta(y^\tau)} \sigma_0(\frac{y^\nu}{\epsilon};\zeta(y^\tau))
  \end{array}
 \right.
\end{equation}
where
\begin{enumerate}[label=(1\alph*)]
 \item \label{theta bullet} $\theta$ is a function of $y^\tau$ only
 
 \item \label{winding number density bullet} $\zeta(y^\tau) := \gamma(\nabla_\tau \theta,\nabla_\tau \theta)$, where $\nabla_\tau$ denotes the tangential gradient along $\Gamma$ and $\gamma_{ij}$ is the induced metric on $\Gamma$ (the ambient metric for this problem is the Minkowski metric - denoted $\eta$).
 
 \item \label{profile bullet} For each $\rho\in\mathbb{R}$ we have that $\Phi_0(\cdot;\rho) := (\phi_0,\sigma_0)(\cdot;\rho):\mathbb{R}\rightarrow\mathbb{R}^2$ satisfies the minimization problem
\begin{gather*}
 \mu(\rho) = \inf\limits_{(f,s)\in\mathcal{A}} \int \left\{\frac{1}{2} \left|(f',s')\right|^2 + V(f,s) + \frac{1}{2}\rho s^2 \right\} \\
 \mathcal{A} := \left\{ (f,s)\in C^1(\mathbb{R},\mathbb{R}^2) \; : \; \lim\limits_{y^\nu\rightarrow \pm\infty} f(y^\nu) = \pm 1 \right\}
\end{gather*} 
 for suitable potentials $V$. In particular, the the profiles $\phi_0$ and $\sigma_0$ in (\ref{formal asymp intro}) are determined by $\zeta(y^\tau)$. 
 
 \item \label{Gamma and Theta EoM bullet}
 $\theta$ and $\Gamma$ satisfy the highly nonlinear, coupled system of PDEs
 \begin{gather}
   \Box_\Gamma \theta = -\gamma\left(\nabla_\tau \log\left[\mu'(\zeta)\right],\nabla_\tau\theta\right)\label{theta intro} \\
   \text{Mean Curvature of }\Gamma = 2 \frac{\mu'(\zeta)}{\mu(\zeta)} \mathrm{I\!I}(\nabla_\tau \theta,\nabla_\tau\theta) \label{Gamma intro}
 \end{gather}
 where the ambient metric that the mean curvature and the second fundamental form $\mathrm{I\!I}$ are defined with respect to is the Minkowski metric.
\end{enumerate}
In this paper, we verify, subject to a non-degeneracy condition, that there does indeed exist a solution to (\ref{equations of motion of interface with a current}) satisfying (\ref{formal asymp intro}) when $n=2$ and when $\Phi$ is an equivariant map.
\bigskip

\begin{figure}[H]
\centering
 \subfloat[Nodal Set]{\includegraphics[scale=0.35]{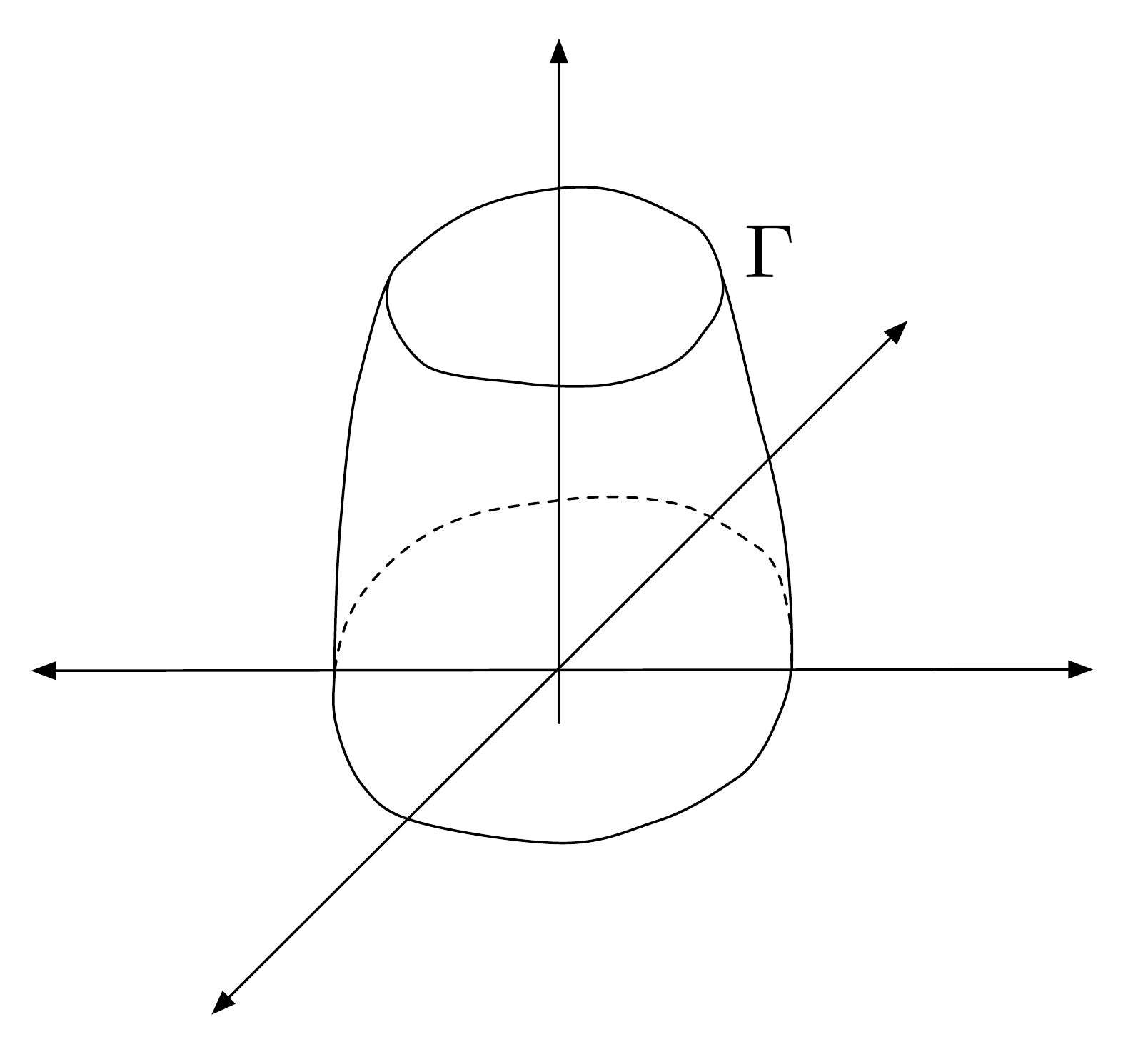}\label{nodal set}}
 \qquad
 \subfloat[Profiles]{\includegraphics[scale=0.65]{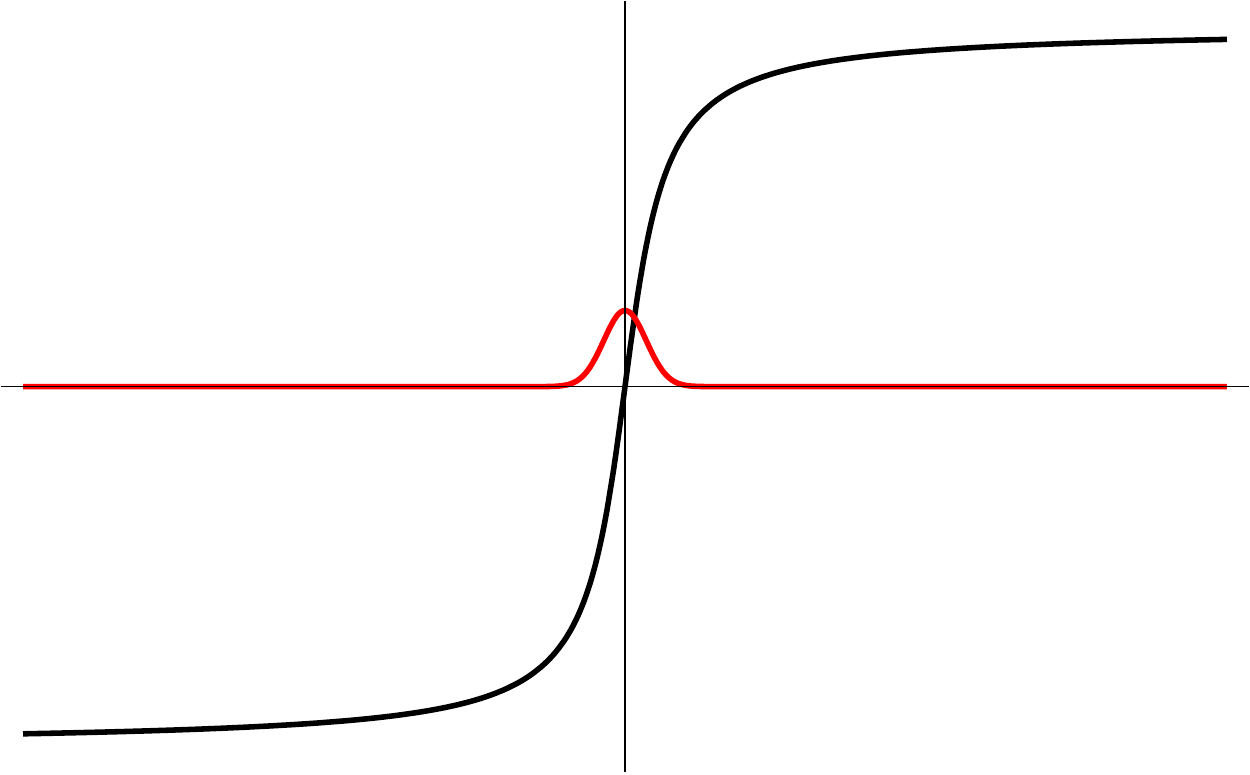}\label{profile image}}
 \caption{The formal asymptotic expansion suggests that there exists a solution $\Phi = (\phi,\sigma)$ to (\ref{equations of motion of interface with a current}) so that for $\theta$ and $\Gamma$ satisfying (\ref{theta intro} - \ref{Gamma intro}), then at each $p\in\Gamma$ we expect that as we move away from $\Gamma$ in the transverse direction $\phi$ looks like the black curve in (\ref{profile image}) and $\sigma$ looks like $e^{\frac{i}{\epsilon}\theta(s_M(p))} \sigma_0$ where $\sigma_0$ looks like the red curve in (\ref{profile image}). Looking at figure \ref{nodal set}, this means that $\sigma$ is exponentially small except near $\Gamma$, $\phi\approx -1$ inside $\Gamma$, $\phi\approx 1$ outside of $\Gamma$, and $\phi$ transitions from $-1$ to $1$ near $\Gamma$.}
\end{figure}

It can be shown that if the winding number density $\gamma^{ij}\partial_i \theta \partial_j \theta$ is sufficiently large, then the $\sigma_0$-field of the approximate solution is $0$. It is believed that there are regimes where a solution may initially have a non-zero current (i.e. $\sigma_0(\frac{\cdot}{\epsilon};\gamma^{ij}\partial_i\theta \partial_j\theta) \neq 0$), but as the system evolves the solution may lose its current. This type of phenomena is referred to as \textbf{current quenching} \cite{vilenkin2000cosmic} and we show in section \ref{current quenching} below that given suitable initial conditions that the solutions we find undergo current quenching.
\bigskip

To the best of our knowledge this is the first paper to consider interface type solutions to a two-component, hyperbolic system.

\subsubsection{Mathematical Background}
There is an extensive mathematical literature with results that are of the type we show in this paper. The unifying theme of these types of results is
\begin{itemize}
 \item There exist solutions to some PDEs which have interfaces, point vortices, or vortex filaments whose dynamics are approximately described by some associated geometric problem.
\end{itemize}
See \cite{jerrard2011semilinear} for a detailed account of these types of results for the scalar elliptic, parabolic, and hyperbolic counterparts of (\ref{equations of motion of interface with a current}). 
\bigskip

Two-component systems have been considered in the physics literature as models for interfaces, point vortices, or vortex filaments in various physical systems \cite{knigavko1998spontaneous,isoshima2002axisymmetric}. However, rigorous mathematical descriptions of solutions to two-component systems of the type we consider are sparse in the math literature. For example, progress on the existence and classification of solutions with interfaces or vortices has been made for various two-component, elliptic systems \cite{alama2006fractional,alama2009structure,alama2012compound,alama2015domain} and (potentially very complicated) ground states of other two-component models subject to physically relevant forcing has been studied \cite{mason2011classification,aftalion2012vortex,aftalion2015thomas-fermi}.
\bigskip

A scalar analogue of (\ref{equations of motion of interface with a current}) is
\begin{equation}\label{scalar hyperbolic case}
  \partial_{tt} u - \Delta u + \frac{1}{\epsilon^2} V'(u) = 0
\end{equation}
where $u:\mathbb{R}^{1+n}\rightarrow \mathbb{R}$ and $V'(u)$ is qualitatively similar to
\[
 \lambda (\left|u\right|^2 - 1)u
\]
In this case, it has been shown that there exists a solutions to (\ref{scalar hyperbolic case}) that have an interface near a codimension one time-like minimal surface \cite{jerrard2011semilinear,galvao2015accelerating}. These results are obtained using weighted energy estimates to show that if one starts with appropriate initial data, then there exists an exact solution to (\ref{scalar hyperbolic case}) that is close to an approximate solution obtained using formal arguments.
\bigskip

One can also consider a version of (\ref{scalar hyperbolic case}) for which $u:\mathbb{R}^{1+n}\rightarrow\mathbb{C}$. In this case, the goal is to find and describe solutions to (\ref{scalar hyperbolic case}) with vortices or vortex filaments. Results describing point vortices and/or vortex filaments in (\ref{scalar hyperbolic case}) and a gauged version of (\ref{scalar hyperbolic case}) have been obtained in \cite{jerrard1999vortex,lin1999vortex,jerrard2011semilinear} and \cite{gustafson2006effective,czubak2015topological}, respectively. Similarly for us, we could consider the case when $\phi:\mathbb{R}^{1+n}\rightarrow \mathbb{C}$. In this case, we would like to find solutions to (\ref{equations of motion of interface with a current}) so that $\phi$ has a vortex filament and $\sigma$ is exponentially small except near the vortex filament of $\phi$, but for now we focus our attention on the case where $\phi:\mathbb{R}^{1+n}\rightarrow\mathbb{R}$ and has an interface.
\bigskip

In contrast to \cite{jerrard2011semilinear,galvao2015accelerating} who use weighted energy estimates, as in \cite{stuart2004geodesic,stuart2004geodesics,gustafson2006effective} we linearize (\ref{equations of motion of interface with a current}) about an approximate solution obtained using a formal asymptotic expansion and we use spectral properties of the linearized operator to show that there exists an exact solution of (\ref{equations of motion of interface with a current}) which is close to the approximate solution. The reason we use a different approach is that in order to resolve the new complexities introduced by the coupling of the current to the interface of $\phi$, a more detailed description of solutions is required that seems hard to obtain using weighted energy estimates.

\subsection{Description of Results}
We will simplify (\ref{equations of motion of interface with a current}) by considering the case when $\Phi = (\phi,\sigma):\mathbb{R}^{1+2}\rightarrow\mathbb{R}\times\mathbb{C}$ is equivariant. That is, we assume $\phi$ and $\sigma$ are of the form
\begin{equation}\label{equivariant assumption}
 \begin{array}{rcl} \phi(t,x) &=& \tilde{\phi}(t,\left|x\right|) \\ \sigma(t,x) &=& e^{i \frac{d}{\epsilon} \operatorname{arg}(x)} \tilde{\sigma}(t,\left|x\right|) \end{array}
\end{equation}
for $(\tilde{\phi},\tilde{\sigma}):\mathbb{R}^{1+1}\rightarrow \mathbb{R}^2$ and $d\in\mathbb{R}/2\pi\epsilon\mathbb{Z}$ is a fixed constant. Using (\ref{equivariant assumption}) to simplify, then for $(t,r)\in \mathbb{R} \times\mathbb{R}_+$ we have that (\ref{equations of motion of interface with a current}) becomes
\[
  \partial_{tt} \left(\begin{array}{c} \tilde{\phi} \\ \tilde{\sigma} \end{array}\right) 
  - \partial_{rr} \left(\begin{array}{c} \tilde{\phi} \\ \tilde{\sigma} \end{array}\right) 
  - \frac{1}{r} \partial_r \left(\begin{array}{c} \tilde{\phi} \\ \tilde{\sigma}\end{array} \right) 
  + \frac{1}{\epsilon^2} \left(\begin{array}{cc} \lambda_\phi (\tilde{\phi}^2 - 1)\tilde{\phi} + \beta \tilde{\sigma}^2 \tilde{\phi} \\ \lambda_\sigma (\tilde{\sigma}^2 - 1)\tilde{\sigma} + \beta \tilde{\phi}^2 \tilde{\sigma}  \end{array}\right) 
  + \frac{1}{\epsilon^2} \left(\begin{array}{cc} 0 & 0 \\ 0 & \frac{d^2}{r^2} \end{array}\right)\left(\begin{array}{c} \tilde{\phi} \\ \tilde{\sigma}\end{array} \right) = \left(\begin{array}{c} 0 \\ 0 \end{array}\right)
\]
We will actually consider the more general family of equations
\begin{equation}\label{equations of motion in polar coordinates}
 \partial_{tt} \left(\begin{array}{c} \tilde{\phi} \\ \tilde{\sigma} \end{array}\right) - \partial_{rr} \left(\begin{array}{c} \tilde{\phi} \\ \tilde{\sigma} \end{array}\right) - \frac{1}{r} \partial_r \left(\begin{array}{c} \tilde{\phi} \\ \tilde{\sigma}\end{array} \right) + \frac{1}{\epsilon^2} \nabla_\Phi V(\tilde{\phi},\tilde{\sigma}) + \frac{1}{\epsilon^2} \left(\begin{array}{cc} 0 & 0 \\ 0 & \frac{d^2}{r^2} \end{array}\right)\left(\begin{array}{c} \tilde{\phi} \\ \tilde{\sigma}\end{array} \right) = \left(\begin{array}{c} 0 \\ 0 \end{array}\right)
\end{equation}
where $V:\mathbb{R}^2\rightarrow\mathbb{R}$ and $\nabla_\Phi V(\phi,\sigma) := (\partial_\phi V, \partial_\sigma V)$. The initial data of (\ref{equations of motion in polar coordinates}) we consider is described in section \ref{Initial Data and the Existence of Solutions} below. We will be interested in solutions to (\ref{equations of motion in polar coordinates}) for the rest of the paper. Hence, we will drop the $\sim$'s from $\tilde{\phi}$ and $\tilde{\sigma}$ for notational convenience.
\bigskip

Consider potentials $V$ satisfying the following assumptions
\begin{equation}\label{potential assumptions}
\begin{array}{ll}
 \text{1.} &  V\in C^2(\mathbb{R}^2,\mathbb{R}) \text{ and } V(\phi,\sigma) = V(\left|\phi\right|,\left|\sigma\right|) \\
 & \\
 \text{2.} & V(\pm 1,0) = 0 \text{ and } V(x,y) > 0 \text{ for all } (x,y)\neq (\pm 1,0)\text{. Furthermore, } V(1,y) < V(x,y) \\
 & \text{ for } x > 1 \text{, } V(x,1) < V(x,y) \text{ for } y > 1\text{, and } V \text{ has a local maximum at } (0,0) \text{ with }  \\ 
 & \nabla_\Phi V\neq 0 \text{ on } (-1,1)\times(-1,1) \text{ otherwise.}  \\
 & \\
 \text{3.} & \left| \operatorname{Hess}_\Phi V(\Phi)\right| \lesssim 1 + \left| \Phi\right|^2 \text{ and }\operatorname{Hess}_\Phi V(\pm 1,0) \geq \lambda_* I \text{ where } I \text{ is the } 2\times 2 \text{ identity } \\
 &  \text{ matrix and } \lambda_* > 0\text{.} \\
 & \\
 \text{4.} & V \text{ satisfies a non-degeneracy and a continuity condition - see (\ref{nondeg con}) below.}
\end{array}
\end{equation}
For potentials $V$ satisfying (\ref{potential assumptions}), we will construct solutions to (\ref{equations of motion in polar coordinates}) so that the $\phi$-field has an interface near a codimension one time-like surface satisfying some geometric problem.
\bigskip

Define 
\begin{equation}\label{shifted potential}
 W(\Phi,R) := V(\Phi) + \frac{1}{2}\frac{d^2}{R^2} \sigma^2
\end{equation}
and we call $W$ the \textbf{shifted potential}. We will denote the gradient of the shifted potential as
\begin{equation}\label{shifted potential gradient}
 w(\Phi,r) := \nabla_\Phi W(\Phi,r)
\end{equation}

Let $\Gamma$ be a codimension one surface parameterized by $(\tau,R(\tau))$ representing the interface of $\phi$ to be determined by (\ref{equations of motion in polar coordinates}).
\begin{lemma}\label{minkowski-normal coordinates existence lemma}
 Let $$\Gamma_T = \left\{ (\tau,R(\tau)) \; : \; 0 \leq \tau \leq T \text{ where } T \text{ is the time of existence of $R$ with } \left| R'\right| < 1\right\}$$ Then there exists a neighbourhood $\mathcal{N}$ (independent of $\epsilon$) of $\Gamma_T$ on which there exists a differentiable solution to
 \begin{equation}\label{minkowksi-normal coordinate existence}
  \left\{
   \begin{array}{cc}
    -\partial_t d_M^2 + \partial_r d_M^2 = 1 & \text{ on } \mathcal{N} \\
    d_M = 0 & \text{ on } \Gamma_T
   \end{array}
   \right.
 \end{equation}
 Furthermore, there exists $s_M:\mathbb{R}^{1+1}\rightarrow \mathbb{R}$ satisfying 
 \begin{equation}\label{minkowski projection}
 \begin{array}{cc}
  -\partial_t d_M \partial_t s_M + \partial_r d_M \partial_r s_M = 0 & \text{ on } \mathcal{N}\\
  (s_M(t,r),R(s_M(t,r))) = (t,r) & \text{ on } \Gamma_T
  \end{array}
 \end{equation}
 so that 
 \begin{equation}\label{characteristic}
  (t,r) = (s_M(t,r),R(s_M(t,r))) + \frac{d_M(t,r)}{\sqrt{1 - R'(s_M(t,r))^2}} (R'(s_M(t,r)),1)
 \end{equation}
\end{lemma}
\noindent A proof for lemma \ref{minkowski-normal coordinates existence lemma} can be found in \cite{jerrard2015dynamics}.
\bigskip

Since $d_M$ is a continuous, then there exists $c > 0$ so that 
\begin{equation}\label{region of main theorem}
 \Sigma_{c,T} := \left\{ (t,r) \; : \; \text{for } 0\leq t \leq T \text{ and } d_M(t,r)\leq c\right\} \subset \mathcal{N}
\end{equation}
where we possibly take $T$ smaller. The initial data of (\ref{equations of motion in polar coordinates}), specified in section \ref{Initial Data and the Existence of Solutions}, is chosen so that $\Phi$ transitions from $(-1,0)$ to $(1,0)$ on $\Sigma_{c,T}$ and so that $\Phi$ is either $(-1,0)$ or $(1,0)$ outside of $\Sigma_{c,T}$.
\bigskip

We look for solutions of (\ref{equations of motion in polar coordinates}) of the form 
\begin{equation}\label{approx soln intro}
 \Phi(t,r) \approx F_0(\frac{d_M
}{\epsilon};R(s_M)) + \epsilon F_1(\frac{d_M}{\epsilon};R(s_M),R'(s_M))
\end{equation}
We could use the same notation as we use in appendix \ref{formal asymptotics} and write $F_0 = F_0(x;\frac{d^2}{R^2})$, but we write $F_0 = F_0(x;R)$ for convenience. As for the $F_1$ term, we tried to show that there exists a solution to (\ref{equations of motion in polar coordinates}) so that $$\Phi \approx F_0(\frac{d_M}{\epsilon};R(s_M))$$ but when $\sigma\neq 0$ the coupling between the $\phi$-field and $\sigma$-field introduces new subtleties into the nature of the solutions that necessitates a more detailed description. Hence, we consider the leading order correction $F_1$. We will see momentarily why $F_1$ depends additionally upon $R'$.  
\bigskip

Plugging $F_0 + \epsilon F_1$ and $r = R(s_M(t,r)) + \frac{d_M(t,r)}{\sqrt{1 - R'(s_M(t,r))^2}}$  into (\ref{equations of motion in polar coordinates}) we find that
\begin{align}
  \frac{1}{\epsilon^2} \text{ term:} & \;\;\; F_0''(\partial_t d_M^2 - \partial_r d_M^2) + \nabla_\Phi W(F_0,R) \label{epsilon^-2 approx solution} \\
  \frac{1}{\epsilon} \text{ term:} & \;\;\; F_1''(\partial_t d_M^2 - \partial_r d_M^2) + \operatorname{Hess}_\Phi W (F_0,R) F_1 + F_0'(\partial_{tt} d_M - \partial_{rr} d_M - \frac{1}{r} \partial_r d_M) \label{epsilon^-1 approx solution} \\
  & - \frac{2}{\sqrt{1 - R'}}\frac{d_M(t,r)}{\epsilon} \left(\begin{array}{cc} 0 & 0 \\ 0 & \frac{d^2}{R^3} \end{array}\right) F_0 + 2 \partial_R F_0' R' (\partial_t d_M \partial_t s_M - \partial_r d_M \partial_r s_M) \nonumber
\end{align}
where $F_0 = F_0(\frac{d_M}{\epsilon};R(s_M))$, $F_1 = F_1(\frac{d_M}{\epsilon};R(s_M),R'(s_M))$, $F_i'$ is the derivative of $F_i$ with respect to the first coordinate $x$, and $\partial_R F_0$ is the derivative of $F_0$ with respect to the second coordinate of $F_0$. There are lower order terms, but these are the two dominate terms. Using the fact that $-\partial_t d_M^2 + \partial_r d_M^2 = 1$, $-\partial_t d_M \partial_t s_M + \partial_r d_M \partial_r s_M = 0$, and the fact that $H(R):= -\partial_{tt} d_M + \partial_{rr} d_M + \frac{1}{r} \partial_r d_M$ is the mean curvature of the surface of rotation generated by $R$ in $\mathbb{R}^{1+2}$, then (\ref{epsilon^-2 approx solution}) and (\ref{epsilon^-1 approx solution}) can be re-written as
\begin{align*}
  \frac{1}{\epsilon^2} \text{ term:} & \;\;\; -F_0'' + \nabla_\Phi W(F_0,R)\\
  \frac{1}{\epsilon} \text{ term:} & \;\;\; -F_1'' + \operatorname{Hess}_\Phi W (F_0,R) F_1 - H(R) F_0' - \frac{2}{\sqrt{1 - R'}}\frac{d_M(t,r)}{\epsilon} \left(\begin{array}{cc} 0 & 0 \\ 0 & \frac{d^2}{R^3} \end{array}\right) F_0
\end{align*}
Heuristically, we find that for $R\in\mathbb{R}$ if $F_0 = (f_0,s_0)(x;R)$ solves
\begin{gather}
 -F_0'' + \nabla_\Phi W(F_0,R) = 0 \label{F_0} \\
 \lim\limits_{x\rightarrow \pm\infty} F_0(x;R) = \pm 1 \nonumber
\end{gather}
and for $L_1(F_0;R)$ defined in (\ref{lin op}) below if $F_1 = (f_1,s_1)(x;R,R')$ solves
\begin{gather}
 L_1(F_0(x;R),R) F_1 = H(R) F_0'(x;R) + 2 \frac{x}{\sqrt{1 - R'^2}} \frac{d^2}{R^3} \left(\begin{array}{c} 0 \\ s_0(x;R) \end{array}\right) \label{F_1} \\
 \lim\limits_{x\rightarrow \pm \infty} F_1(x;R,R') = 0 \nonumber
\end{gather}
then $F_0(\frac{d_M}{\epsilon};R(s_M)) + \epsilon F_1(\frac{d_M}{\epsilon};R(s_M),R'(s_M))$ has the properties that we are looking for and looks to be a good approximate solution. Note that $F_0$ and $F_1$ depend on this so far unknown function $R$ parameterizing $\Gamma$. In fact, $F_0$ depends on $R$ and $F_1$ depends on $R$, $R'$, and $R''$. We will see momentarily that $F_1$ actually only depends on $R$ and $R'$.
\bigskip

Differentiating (\ref{F_0}) with respect to $x$ yields
\[
 - (F_0')'' + \operatorname{Hess}_\Phi W(F_0;R) F_0' = 0
\]
Define
\begin{equation}\label{lin op}
  L_1 := -\frac{d^2}{dx^2} I_{2\times 2} + \operatorname{Hess}_\Phi W(F_0;R)
\end{equation}
We assume that 
\begin{equation}\label{nondeg con}
 \ker(L_1(F_0;R)) = \operatorname{span}\left\{ F_0'\right\}
\end{equation}
holds for $R\in (r_0,r_1)$ with $0 \leq r_0 < r_1 \leq \infty$. This is the non-degeneracy condition we assume $V$ satisfies. Further, we also assume that for $F_0(\cdot;R)$ satisfying (\ref{F_0}),
\begin{gather}\label{cont con}
 R\mapsto F_0(\cdot;R)
\end{gather}
is a continuous map for all $R\in (r_0,r_1)$. This is the continuity condition we assume $V$ satisfies. Of particular interest to us are potentials $V$ for which (\ref{nondeg con}) holds for a range of $R$ for which $s_0(\cdot;R)\neq 0$.
\bigskip

The reason we assume (\ref{nondeg con}) and (\ref{cont con}) is they allow us to conclude that $F_0(\cdot;R)$ is actually $C^2$ in $R$. They also give us certain spectral estimates (see theorem \ref{spectral estimate theorem}) which will be important in the proof of theorem \ref{main theorem intro}. Further, we will need to use (\ref{nondeg con}) to find $F_1$ solving (\ref{F_1}) and decaying at infinity in proposition \ref{existence of F_1}. In fact, to show that such an $F_1$ exists, we need to solve an equation like $L_1(F_0;R) f = g$ with $g\in L^2(\mathbb{R};\mathbb{R}^2)$. A necessary condition for this to be solvable is that $g\in \ker(L_1(F_0;R))^\perp$ where $\perp = \perp_{L^2}$. This necessary condition plus (\ref{nondeg con}) suggests to us that for (\ref{approx soln intro}) to hold, then $R$ must solve 
\begin{gather}
 H(R) \int \left| F_0'(x;R) \right|^2 dx - \frac{1}{\sqrt{1 - R'^2}} \frac{d^2}{R^3} \int \left| s_0(x;R)\right|^2 dx = 0 \label{R} \\
 R(0) \in (r_0,r_1)  \text{ and }  R'(0) = 0 \nonumber
\end{gather}
Since $R(s_M)$ could leave $(r_0,r_1)$, we have that the approximate solution $F_0(\frac{d_M}{\epsilon};R(s_M)) + \epsilon F_1(\frac{d_M}{\epsilon};R(s_M),R'(s_M))$ is only valid up to some finite time $T$ as $F_1$ is only guaranteed to exist for as long as $R(s_M)\in (r_0,r_1)$. Furthermore, the mean curvature $H(R)$ of $\Gamma$ contains a $R''$ term. Using (\ref{R}), one can express $R''$ in terms of $R$ and $R'$. From this we see that $F_1$ actually only depends on $R$ and $R'$ as stated earlier. 
\bigskip

The main result obtained in this paper shows that $F_0(\frac{d_M}{\epsilon};R) + \epsilon F_1(\frac{d_M}{\epsilon};R,R')$ is indeed a good approximate solution. 
\begin{theorem}\label{main theorem intro}
 Given suitable initial conditions (see section \ref{Initial Data and the Existence of Solutions} below), there exists a solution $\Phi$ to (\ref{equations of motion in polar coordinates}), a function $a:\mathbb{R}\rightarrow \mathbb{R}$, and constants $\overline{T} \leq T$ and $\delta > 0$ (both independent of $\epsilon$) with $$[0,\overline{T}]\times [\mathbb{R}_+ \setminus \left(R(0) -\delta,R(0)+\delta\right)] \cup \Sigma_{c,\overline{T}} = [0,\overline{T}]\times \mathbb{R}_+$$ so that  on $\Sigma_{c,\overline{T}}$, $R(s_M(t,r))\in (r_0,r_1)$ and
 \begin{equation}\label{main estimate on Sigma space}
  \left\| \Phi - F_0(\frac{d_M - a(s_M)}{\epsilon};R(s_M)) - \epsilon F_1(\frac{d_M - a(s_M)}{\epsilon};R(s_M),R'(s_M)) \right\|_{L_t^1 H_r^1(\Sigma_{c,\overline{T}})} \lesssim \epsilon^2
 \end{equation}
 \begin{equation}\label{main estimate on Sigma time}
  \left\| \partial_t \left[\Phi - F_0(\frac{d_M - a(s_M)}{\epsilon};R(s_M)) - \epsilon F_1(\frac{d_M - a(s_M)}{\epsilon};R(s_M),R'(s_M))\right] \right\|_{L_t^1 L_r^2(\Sigma_{c,\overline{T}})} \lesssim \epsilon^2
 \end{equation}
 where $d_M$ and $s_M$ are functions of $(t,r)$ and off of $\Sigma_{c,\overline{T}}$
 \begin{eqnarray*}
  \Phi = (-1,0) & \text{ for } & (t,r)\in [0,\overline{T}]\times \left[0,R(0) - \delta\right) \\
  \Phi = (1,0) & \text{ for } & (t,r)\in [0,\overline{T}]\times \left(R(0) - \delta,\infty\right)
 \end{eqnarray*}
\end{theorem}

Most importantly, this theorem tells us that there exists a solution $\Phi$ to (\ref{equations of motion in polar coordinates}) with the properties that we want. Namely, there exists a solution $\Phi = (\phi,\sigma)$ so that $\phi$ has an interface and for appropriate potentials $V$, $\sigma$ is exponentially small except near the interface of $\phi$.

\subsubsection{Initial Data and the Existence of Solutions}\label{Initial Data and the Existence of Solutions}
Provided suitable initial data, to be described shortly, showing that the superconducting interface model is globally well posed is a standard exercise as this problem is energy subcritical \cite{shatah1998geometric,tao2006nonlinear}. The initial data we consider is described next.
\bigskip

It turns out that the Minkowskian distance $d_M$ from lemma \ref{minkowski-normal coordinates existence lemma} is not necessarily defined everywhere. It is, however, defined on the set $\Sigma_{c,T}$, defined in (\ref{region of main theorem}), for some $c$, $T > 0$. For $0 < \overline{T} \leq T$, define $$\mathcal{B} := [0,\overline{T}]\times (R(0) - \delta,R(0) + \delta)$$ where $\delta > 0$ is chosen so that $\mathcal{B} \subset \Sigma_{c,T}$. We choose the initial data of $\Phi$ as
\begin{align}
 \Phi(0,r) = \left\{\begin{array}{cl} (-1,0) & \text{ for } 0 \leq r < b_1 \\ (1,0) & \text{ for } r > b_2 \end{array}\right. \label{initial data of Phi} \\
 \partial_t \Phi(0,r) = 0 \text{ for } 0 \leq r < b_1 \text{ and } r > b_2 \label{initial data of Phi_t}
\end{align}
for some $0 < b_1 < b_2$ to be chosen shortly. Since (\ref{equations of motion in polar coordinates}) is a wave equation, there is a finite speed of propagation of data. We choose $b_1$ and $b_2$ so that
\[
 \Phi(t,r) = \left\{\begin{array}{cl} (-1,0) & \text{ for } 0 \leq r < R(0) - \delta \\ (1,0) & \text{ for } r > R(0) + \delta \end{array}\right.
\]
for all $0 \leq t \leq \overline{T}$. Thus, we know the value of $\Phi$ outside of $\Sigma_{c,\overline{T}}$. This reduces the analysis to controlling the error between $\Phi$ and the right hand side of (\ref{approx soln intro}) on the region where $\Phi$ transitions from $(-1,0)$ to $(1,0)$ - the region $\Sigma_{c,\overline{T}}$.

\begin{figure}[H]
 \centering
 \includegraphics[scale = 0.5]{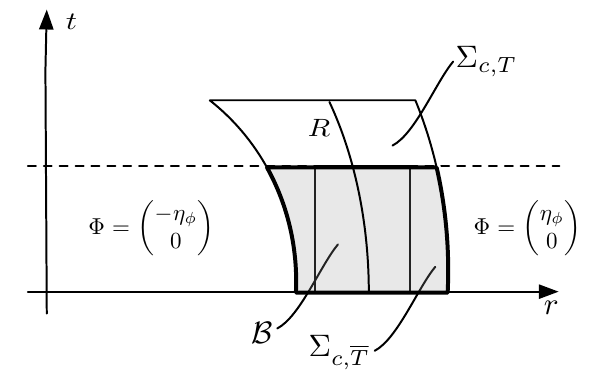}
 \caption[font={small}]{$\Sigma_{c,T}$ is a neighbourhood on which the Minkowski normal coordinates are well defined. Pick a rectangle $\mathcal{B} = [0,\overline{T}]\times (R(0) - \delta,R(0) + \delta)$ that lays within $\Sigma_{c,T}$. We then pick the initial data of (\ref{equations of motion in polar coordinates}) so that $\Phi = (-1,0)$ to the left of $\mathcal{B}$ and $\Phi = (1,0)$ to the right of $\mathcal{B}$.}\label{figure of initial data}
\end{figure}

\section{Physical Motivation: Superconducting Strings}\label{physics background}
Motivated by \cite{nielsen1973vortex}, Witten introduced a two-component model, closely related to the abelian-Higgs model, to describe finite energy solutions with vortex filaments supporting superconducting currents \cite{witten1985superconducting}. We call this model the \textbf{superconducting string model}. It was our initial consideration of this model that lead us to study (\ref{equations of motion of interface with a current}) - the superconducting interface model. We will describe what lead us to consider (\ref{equations of motion of interface with a current}), but in order to do so we will first need to describe the superconducting interface model.
\bigskip

In \cite{witten1985superconducting}, an effective action for the superconducting string model using formal arguments was derived. The effective action found suggest that
\begin{enumerate}[label=(2\alph*)]
 \item \label{filaments} there should be solutions to the model with a vortex filament with a superconducting current
 
 \item \label{filaments dynamics} the vortex filament is near a codimension two time-like surface $\Gamma$, where $\Gamma$ satisfies a geometric equation that is coupled in a highly nonlinear way to the phase of the current and an ambient vector potential representing an external electromagnetic field
\end{enumerate}
To obtain this effective action, it is proposed that there there exists solutions to the superconducting string model whose profiles to leading order only depend on $d_M$. In contrast, the ansatz we use to derive an effective action depends additionally on the gradient of the phase of the field corresponding to the current, see (\ref{formal asymp intro}). The effective action derived in the physics literature in the case when the phase of the current is decoupled from the vector potential looks like the effective action we derived in (\ref{effective action 2}) with $\mu(\gamma^{ij} \partial_i\theta\partial_j\theta)$ replaced with the first order Taylor approximation of $\mu$ about $0$. Results obtained in this paper suggest that, at least for the superconducting interface model, that the physics ansatz leads to a less accurate approximation of solutions.
\bigskip

To illustrate how the superconducting interface model is related to the superconducting string model, we first need to state the superconducting string model. However, before we can state the superconducting string model we need some notation. We will denote the complex scalar fields as $\phi,\sigma:\mathbb{R}^{1+3}\rightarrow \mathbb{C}$ and denote their associated gauge fields as $A_\phi,A_\sigma:\mathbb{R}^{1+3}\rightarrow \mathbb{R}^4$. We define the covariant derivatives associated to the $\phi$ and $\sigma$ fields as $\nabla_\phi = \nabla - i q_\phi A_\phi$ and $\nabla_\sigma = \nabla - i q_\sigma A_\sigma$, respectively, where $q_\phi$, $q_\sigma\in \mathbb{R}$ are the coupling constants between $(\phi,\sigma)$ and their associated gauge fields. As is standard notation, we define $F_{\phi,\mu\nu} := \partial_\mu A_{\phi,\nu} - \partial_\nu A_{\phi,\mu}$ and similarly define $F_{\sigma,\mu\nu}$. Finally, for $(\lambda_\phi,\lambda_\sigma,\beta)\in\mathbb{R}_+^3$ the \textbf{superconducting string 
potential} is
\begin{equation}\label{potential}
 V_S(\phi,\sigma) = \frac{\lambda_\phi}{4}(\left|\phi\right|^2 - 1)^2 + \frac{\lambda_\sigma}{4}(\left|\sigma\right|^2 - 2) \left|\sigma\right|^2 + \frac{\beta}{2}\left|\phi\right|^2 \left|\sigma\right|^2
\end{equation}
\smallskip

The Lagrangian of the superconducting string model is defined as
\begin{equation}\label{witten's model}
 \mathcal{L} = \frac{1}{2}\eta^{\alpha\beta} \overline{\nabla_{\phi,\alpha} \phi} \nabla_{\phi,\beta} \phi + \frac{1}{2}\eta^{\alpha\beta} \overline{\nabla_{\sigma,\alpha} \sigma} \nabla_{\sigma,\beta} \sigma + \frac{1}{\epsilon^2} V_S(\phi,\sigma) + \frac{\epsilon^2}{4} F_{\phi,\mu\nu} F_\phi^{\mu\nu} + \frac{\epsilon^2}{4} F_{\sigma,\mu\nu} F_\sigma^{\mu\nu}
\end{equation}
where $0 < \epsilon \ll 1$ and $\eta = \operatorname{diag}(-1,1,1,1)$ is the Minkowski metric. An important feature of this model that is worth highlighting is that the $\sigma$-field has a $U(1)$ gauge symmetry. See \cite{vilenkin2000cosmic} for an in depth discussion of the physics behind this model.
\bigskip

To obtain (\ref{equations of motion of interface with a current}), two changes to the superconducting string model will be made. The first is to consider $\phi:\mathbb{R}^{1+n}\rightarrow\mathbb{R}$.  In this case, $\phi$ loses its $U(1)$ gauge symmetry and gains a discrete symmetry. In particular, this allows for $\phi$ to have an interface. The second change we make is to simplify the problem by decoupling the current from the ambient vector potential. To do this, set $q_\sigma = 0$. Applying these changes to (\ref{witten's model}), one obtains the Lagrangian for the superconducting interface model
\begin{equation}\label{superconducting interface model}
 \mathcal{L} = \frac{1}{2}\eta^{\alpha\beta} \partial_\alpha \phi \partial_\beta\phi + \frac{1}{2}\eta^{\alpha\beta}\overline{\partial_\alpha \sigma}\partial_\beta\sigma + \frac{1}{\epsilon^2} V_S(\phi,\sigma)
\end{equation}

\section{Effective Equations}
\subsection{Change to Minkowski Normal Coordinates}
Suppose $\Gamma$ is a codimension 1 time-like surface parameterized by $\left\{(y^0,R(y^0))\right\}$ representing the interface of the $\phi$-field we wish to find. Note that (\ref{equations of motion in polar coordinates}) will determine $\Gamma$. It turns out that doing a change of coordinates from polar coordinates to a coordinate system centred about $\Gamma$ to ``straighten out'' $\Gamma$ is quite useful. In fact, it is in these new coordinates that we will find $\Phi$ which has the properties we want. 
\bigskip

The coordinate system we will work in, called \textbf{Minkowski normal coordinates} and denoted $(y^0,y^1)$, are defined as
\begin{equation}\label{m-normal coordinates} 
(t,r) = (y^0,R(y^0)) + y^1 \nu(y^0) 
\end{equation}
where
\[
 \begin{array}{rcl}
  \eta(\partial_{y^0} (y^0,R(y^0)),\nu(y^0)) = 0 & \text{ and } & \eta(\nu(y^0),\nu(y^0)) = 1
 \end{array}
\]
We have a choice of $\nu(y^0)$ and pick  
\[
 \nu(y^0) := \frac{1}{\sqrt{1 - (R')^2}}\left(R', 1\right)
\]
Going back to lemma \ref{minkowski-normal coordinates existence lemma} in the introduction, $y^0$ can be identified with $s_M$ and $y^1$ can be identified with $d_M$.
\bigskip

The action integral associated to (\ref{equations of motion in polar coordinates}) is
\begin{equation}\label{superconducting interface model polar coordinates}
 S(\Phi) := \int \left\{-\frac{1}{2} \partial_t \Phi^2 + \frac{1}{2} \partial_r \Phi^2 + \frac{1}{\epsilon^2} W(\Phi,r)\right\} r dt dr
\end{equation}
where $\Phi = (\phi,\sigma)$ and where $W$ is the shifted potential defined in (\ref{shifted potential}). Define 
\begin{equation}\label{m and n}
 \begin{array}{ccc} m := (1 - (R')^2)^{-1/2} & \textrm{ and } & n := 1 + y^1 m^3 R'' \end{array}
\end{equation}
A computation shows that
\begin{equation}\label{change of coordinates}
  \left(\begin{array}{c} \partial_t \\ \partial_r \end{array}\right)
 = \frac{m}{n} \left(\begin{array}{cc} m & -n\;R' \\ -m\;R' & n \end{array}\right)\left(\begin{array}{c} \partial_{y^0} \\ \partial_{y^1} \end{array}\right)
\end{equation}
In Minkowski normal coordinates, $S$ is
\begin{equation}\label{superconducting interface model in minkowski normal coordinates}
 S(\Phi) = \int \left\{ -\frac{m^2}{2 n^2} \partial_{y^0} \Phi^2 + \frac{1}{2} \partial_{y^1} \Phi^2 + \frac{1}{\epsilon^2} W(\Phi,R(y^0) + y^1 m(y^0))\right\} n(y^0,y^1) (R(y^0) + y^1 m(y^0)) dy^0 dy^1
\end{equation}
The equations of motion of (\ref{superconducting interface model in minkowski normal coordinates}) are
\begin{equation}\label{Wave equation in minkowski-normal coordinates}
 \frac{m^2}{n^2}\partial_{y^0 y^0}\Phi + B^\alpha \partial_\alpha \Phi - \partial_{y^1 y^1}\Phi + \frac{1}{\epsilon^2} w(\Phi, R(y^0) + y^1 m(y^0)) = 0
\end{equation}
with initial data as described in section \ref{Initial Data and the Existence of Solutions} and where $w$ was defined in (\ref{shifted potential gradient}) and we've defined
\begin{equation}\label{B^0}
 B^0 := \frac{m}{n}\partial_{y^0}(\frac{m}{n}) + \frac{1}{(R + y^1 m)}\frac{m^2}{n}R'
\end{equation}
\begin{equation}\label{B^1}
 B^1 := - \frac{m^3}{n}R'' - \frac{1}{(R + y^1 m)} m
\end{equation}
\smallskip

It turns out that Minkowski normal coordinates are not well defined everywhere. They are, however, well defined on $[0,y_*^0]\times [-y_*^1,y_*^1]$ where $y_*^0$ and $y_*^1$ are determined by the time of existence of $R$. We will also choose $y_*^0$ possibly smaller so that $R(y^0)\in (r_0,r_1)$ for $y^0\in [0,y_*^0]$ where $r_0 < r_1$ come from the non-degeneracy condition (\ref{nondeg con}). Using our choice of initial data, see section \ref{Initial Data and the Existence of Solutions}, $y_*^1$ can be chosen possibly smaller so that for $0\leq y^0 \leq y_*^0$, then $\Phi(y^0,y^1) = (-1,0)$ for $y^1 < -y_*^1$ and $\Phi(y^0,y^1) = (1,0)$ for $y^1 > y_*^1$. Thus, we are left to find solutions to (\ref{superconducting interface model in minkowski normal coordinates}) on this neighbourhood connecting these two states. For the rest of the paper we consider (\ref{superconducting interface model in minkowski normal coordinates}).

\subsection{Expansion}
We would like to find a solution to (\ref{Wave equation in minkowski-normal coordinates}) which has an interface which is centred on a function $R(y^0)$. Furthermore, we'd like the solution to only have an $O(1)$ change for $O(\epsilon)$ movements in transverse directions and an $O(1)$ change for $O(1)$ movements in tangential directions.
\bigskip

We would like to construct a solution $\Phi$ of (\ref{Wave equation in minkowski-normal coordinates}) as 
\begin{equation}\label{approx solution w/ leading order correction}
 \Phi \approx F_0(\frac{y^1}{\epsilon};R) + \epsilon F_1(\frac{y^1}{\epsilon};R,R')
\end{equation}
where $F_i = (f_i,s_i)$ and each $F_i$ is independent of $\epsilon$ (as we discussed following (\ref{approx soln intro})). If we plug $\Phi$ into (\ref{Wave equation in minkowski-normal coordinates}) and expand, we can find the equations that each $F_i$ must satisfy. 
\bigskip

Doing so, we find that $F_0 = F_0(y^1;R)$ must satisfy
\begin{equation}\label{F_0 equation}
 -\partial_{y^1 y^1} F_0 + w(F_0,R)=0
\end{equation}
for each $R\in\mathbb{R}$ where $w$ was defined in (\ref{shifted potential gradient}). Since $w(\cdot,R)$ depends on $R$, we can see why solutions $F_0$ depend on $R$ too.
\bigskip

We also find that for $v=R'$, then $F_1 = F_1(y^1;R,v)$ must satisfy
\begin{equation}\label{F_1 equation}
 L_1(F_0,R) F_1 = H(R) \partial_{y^1} F_0 - \frac{y^1}{\sqrt{1 - (R')^2}} \partial_R w(F_0,R)
\end{equation}
where 
\begin{equation}\label{linearized operator}
 L_\epsilon(F_0,R) := -\partial_{y^1 y^1} + \frac{1}{\epsilon^2} \operatorname{Hess}_\Phi W(F_0,R)
\end{equation}
\begin{equation}\label{mean curvature}
 H(R) := \frac{1}{\sqrt{1 - (R')^2}}\left(\frac{R''}{1 - (R')^2} + \frac{1}{R}\right)
\end{equation}
\begin{equation}\label{partial w}
 \partial_R w(\Phi,r) := \lim\limits_{\Delta r\rightarrow 0} \frac{w(\Phi,r+\Delta r) - w(\Phi,r)}{\Delta r}
\end{equation}
The operator $L_\epsilon(F_0,R)$ is the linearized operator of (\ref{superconducting interface model in minkowski normal coordinates}), linearized about $F_0$, and $H(R)$ is the \textbf{mean curvature} of the surface of rotation generated by $R$ in $\mathbb{R}^{1+2}$. A necessary condition for (\ref{F_1 equation}) to be solvable is that the right hand side of (\ref{F_1 equation}) must be orthogonal to the kernel of $L_\epsilon(F_0,R)$. This implies that 
\begin{equation}\label{R equation}
 H(R) \int\limits_\mathbb{R} \partial_{y^1} F_0(\cdot;R)^2 - \frac{d^2}{R^3} m \int\limits_\mathbb{R} s_0(\cdot;R)^2 = 0
\end{equation}
\bigskip

The following is an important estimate regarding the operator $L_\epsilon(F_0;R)$ that will be used to verify that there exists a solution to (\ref{Wave equation in minkowski-normal coordinates}) satisfying (\ref{approx solution w/ leading order correction}).
\begin{theorem}[Spectral Estimate]\label{spectral estimate theorem}
Suppose $F_0$ and $R$ satisfy (\ref{F_0 equation}) and (\ref{R equation}), respectively. By assumption 4 of (\ref{potential assumptions}), $\ker(L_\epsilon(F_0,R)) = \operatorname{span}\left\{ \partial_{y^1} F_0 \right\}$. In particular, this implies that for $\perp = \perp_{L^2}$, then for any $\xi\in \ker\left(L_\epsilon(F_0,R) \right)^\perp$ we have
 \begin{equation}\label{spectral estimate}
  \frac{1}{\epsilon^2} \left\|\xi\right\|_{L^2(\mathbb{R})}^2 \lesssim \int\limits_{\mathbb{R}} \xi\cdot L_\epsilon(F_0;R)\xi
 \end{equation}
\end{theorem}

\noindent \underline{Sketch of Proof}: For fixed $R\in (r_0,r_1)$ define $$X = \left\{ \xi\in H^1(\mathbb{R};\mathbb{R}^2)\; :\; \left\|\xi\right\|_2 = 1 \text{ and } \left<\xi,\partial_{y^1} F_0\right>_2 = 0\right\}$$ $$I(\xi) := \int\limits_{\mathbb{R}} \xi\cdot L_\epsilon(F_0,R)\xi$$ To see why (\ref{spectral estimate}) holds, we want to show that $$m := \inf\limits_{\xi\in X} I(\xi) > 0$$ Clearly, if $m \geq \lambda_*$, where $\lambda_*$ is from assumption 3 of (\ref{potential assumptions}), then there is nothing to be done. Suppose $m < \lambda_*$. If $m=0$ and there exists $\xi\in X$ at which this infimum is attained then by the non-degeneracy condition (\ref{nondeg con}) $\xi\propto \partial_{y^1} F_0$. Since $\xi\in X$, then $\xi\perp \partial_{y^1} F_0$ which implies that $\xi=0$. This contradicts the fact that $\left\| \xi\right\|_2 = 1$. Thus, if we show that there exists $\xi\in X$ at which the infimum of $I$ is attained, then we are done.
\bigskip

Let $\xi_n\in X$ be a minimizing sequence. We have that 
\begin{gather}\label{H^1 uniform bound}
\left\| \xi_n'\right\|_2 \leq I(\xi_n) < C \;\;\; \text{ and } \;\;\; \left\| \xi_n\right\|_2 = 1 \leq C
\end{gather}
Thus, by possibly passing to a subsequence we have that $\xi_{n_k}\rightharpoonup \xi$ in $H^1$. Furthermore, for $\left[f\right]_{1/2}$ denoting the Holder-$1/2$ constant of $f$ we have that $$\left[\xi_{n_k}\right]_{1/2} \leq \left\| \xi_n'\right\|_2 \leq I(\xi_n) < C$$ and so $\xi_{n_k}\rightarrow \xi$ locally uniformly. Thus, we have that $$I(\xi) \leq \lim\inf_{n_k\rightarrow \infty} I(\xi_{n_k})$$
\smallskip

Note that we still have that $\xi\perp\partial_{y^1} F_0$. Suppose $\left\|\xi\right\|_2 = t\in\left[0,1\right]$. We have that $$I(\xi) = t^2 I(\frac{\xi}{t}) \geq t^2 m$$ Using concentration compactness type arguments one can show that $$m = \lim\inf_{n\rightarrow \infty} I(\xi_n) \geq t^2 m + (1 - t^2) \lambda_* \geq m$$ with equality if and only if $t = 1$. Thus, $t = 1$ and so the infimum of $I$ is attained in $X$.
\begin{flushright}
$\Box$
\end{flushright}

\subsection[Existence of F0, F1, and R]{Existence of $F_0$, $F_1$, and $R$}
Set 
\begin{gather}
 \mu(R) = \inf_{(f,s)\in\mathcal{A}} \int\limits_{\mathbb{R}} \mu(f,s;R) \label{minimization problem} \\
 \mathcal{A} = \left\{ (f,s)\in H_{loc}^1(\mathbb{R}) \; : \; \int\limits_{\mathbb{R}} \mu(f,s;R) < \infty,\; f(0) = 0,\; f(\pm \infty) = \pm 1 \right\} \label{min set}
\end{gather}
where
\begin{equation}\label{min arg}
 \mu(f,s;R) := \frac{1}{2}(f')^2 + \frac{1}{2}(s')^2 + W(f,s;R)
\end{equation}
Notice that without the requirement that $f(0) = 0$, then any translation of a minimizer is another minimizer. This condition kills this degeneracy. We will now show that there exists $(f,s)\in\mathcal{A}$ at which $\mu(f,s;R)$ attains its infimum.
\begin{proposition}\label{existence of minimizers}
 There exists $(f,s)\in\mathcal{A}$ that solves the minimization problem (\ref{minimization problem}).
\end{proposition}

Minimization problems like proposition \ref{existence of minimizers} have been studied extensively. Arguments used in \cite{alikakos2008connection} can be used to prove this proposition. In particular, we have the following corollary
\begin{corollary}\label{A tilde corollary}
 Set $$\tilde{\mathcal{A}} = \left\{ (f,s)\in \mathcal{A} \; : \; f(y^1) = 0 \textrm{ iff } y^1 = 0, \left|f\right| \leq 1, 0 \leq s \leq 1 ,f \textrm{ is odd}, s \textrm{ is even} \right\}$$ then there exists a minimizer of $$\inf\limits_{(f,s)\in\tilde{\mathcal{A}}} \mu(f,s;R)$$
\end{corollary}

The proof of this corollary follows from the following lemma
\begin{lemma}\label{minimizer properties}
 \begin{gather*}
  \inf\limits_{(f,s)\in \mathcal{A}} \int\limits_{\mathbb{R}} \mu(f,s;R) = \inf\limits_{(f,s)\in\tilde{\mathcal{A}}} \int\limits_{\mathbb{R}} \mu(f,s;R)
 \end{gather*}
\end{lemma}
The details of this proof are omitted, but the idea is to modify $(f,s)\in\mathcal{A}$ by using appropriate symmetrizations and translations to show that there is $(\tilde{f},\tilde{s})\in\tilde{\mathcal{A}}$ so that $$\int\limits_{\mathbb{R}} \mu(\tilde{f},\tilde{s};R) \leq \int\limits_{\mathbb{R}} \mu(f,s;R)$$ For notational convenience, we drop the $\sim$ from $\tilde{\mathcal{A}}$. 
\bigskip

For each fixed $R$ corollary \ref{A tilde corollary} gives us the existence of a minimizer $F(\cdot;R)$. One can verify using the non-degeneracy condition (\ref{nondeg con}) along with the continuity condition (\ref{cont con}) that $$R\mapsto F(\cdot;R)$$ is actually a $C^2$ map. In fact, for this to be true it suffices that $F(\cdot;R)$ is a local minimizer, modulo discrete symmetries, for $R\in (r_1,r_2)$. Since $R(0)\in (r_0,r_1)$, we can then plug $F(\cdot;R)$ into (\ref{R equation}) and solve for $R$ using standard ODE techniques \cite{coddington1955theorem} as long as $R$ remains in $(r_0,r_1)$.
\bigskip

\begin{proposition}\label{existence of F_1}
 Define
 \begin{equation}\label{g mean curve}
  g(R,v) := \frac{1}{\sqrt{1 - v^2}} \frac{d^2}{R^3} \frac{\left\| s_0(\cdot;R)\right\|_{L^2(\mathbb{R})}^2}{\left\| \partial_{y^1} F_0(\cdot;R)\right\|_{L^2(\mathbb{R})}^2}
 \end{equation}
 For every $(R,v)\in (r_0,r_1)\times (-1,1)$, there exists a unique $F_1 = F_1(y^1;R,v)\in H^1(\mathbb{R};\mathbb{R}^2)$ solving
 \begin{gather*}
  L_1(F_0,R) F_1 = g(R,v) \partial_{y^1} F_0(\cdot;R) - \frac{y^1}{\sqrt{1-v^2}} \partial_R w(F_0,R) \\
  \int F_1 \cdot \partial_{y^1} F_0(\cdot;R) = 0
 \end{gather*}
 where $L_1(F_0,R)$ was defined in (\ref{linearized operator}) and $\partial_R w(F_0,R)$ was defined in (\ref{partial w}).
\end{proposition}

\noindent \underline{Proof of Proposition \ref{existence of F_1}}: The main tool in proving the existence of $F_1$ is the spectral theorem \cite{reed1980functional} applied to the unbounded operator $L_1 (F_0,R):L^2(\mathbb{R};\mathbb{R}^2)\rightarrow L^2(\mathbb{R};\mathbb{R}^2)$. To use the spectral theorem, we first need to describe the spectrum of $L_1(F_0,R)$. The essential spectrum of $L_1(F_0,R)$ is $$\sigma_{ess}(L_1(F_0,R)) = \left[\lambda_*,\infty\right)$$ where $\lambda_* > 0$ is from assumption 3 of (\ref{potential assumptions}). Since $L_1(F_0,R)$ is self-adjoint, then by the spectral theorem there exists a spectral projection $E_\lambda$ of $L_1(F_0,R)$. Since there exists a spectral gap, by theorem \ref{spectral estimate theorem}, and since $$g(R,v) \partial_{y^1} F_0(\cdot;R) + \frac{y^1}{\sqrt{1-v^2}} \partial_R w(F_0,R) \perp \ker(L_1(F_0,R))$$ then $F_1\in H^1$ satisfying
\begin{equation}\label{spec F1 existence}
 \left<\psi,F_1\right> := \int\limits_{\alpha_0}^\infty \frac{1}{\lambda} d\left<\psi, E_\lambda\left(g(R,v) \partial_{y^1} F_0(\cdot;R) + \frac{y^1}{\sqrt{1-v^2}} \partial_R w(F_0,R)\right)\right>
\end{equation}
for all $\psi\in H^1$ solves (\ref{F_1 equation}), where $0 < \alpha_0$ is the second smallest value of the spectrum of $L_1(F_0,R)$. Moreover, from (\ref{spec F1 existence}), one has that $F_1\perp \partial_{y^0} F_0$.
\begin{flushright}
 $\Box$
\end{flushright}

\subsection{Properties of Profiles}\label{properties}
\subsubsection[Regimes where s0 is nonzero]{Regimes where $s_0$ is nonzero}
An important feature of what we are studying is that our model is an interface with a current. For a current to exist, it is necessary that $s_0 \neq 0$. For this section we fix the potential to be (\ref{potential}). It can easily be checked that for $\lambda_\sigma < \beta < \lambda_\phi$, then this potential satisfies all the conditions set out in (\ref{potential assumptions}), except for condition 4. We believe that (\ref{potential}) satisfies this as well, but we did not verify this. We will show that if $(f_0,s_0)$ satisfies the minimization problem set out in (\ref{minimization problem}) for this potential, then $s_0\neq 0$.
\bigskip

Set $$E(f,s) = \int\limits_{\mathbb{R}} \left\{ \frac{1}{2}(f')^2 + \frac{1}{2} (s')^2 + W(f,s,R)\right\} dy^1$$ For $s=0$, then it is known that $E_0(f) := E(f,0)$ is uniquely minimized when  $$f_{min}(y^1) = \tanh(\sqrt{\frac{\lambda_\phi}{2}} y^1)$$ The goal is to find parameters so that $E(f_{min},s) < E(f_{min},0)$ for some non-zero $s$ with $(f_{min},s)\in \mathcal{A}$.
\bigskip

\begin{proposition}\label{non-zero s minimizer}
 If the constants of the model additionally satisfy
 \begin{equation}\label{non-zero s minimizer further constraint}
  \beta < \frac{3}{2} \lambda_\sigma
 \end{equation}
 then for sufficiently large $R$ there exists a minimizer $(f,s)\in\mathcal{A}$ of (\ref{minimization problem}) with $s\neq 0$.
\end{proposition}

\noindent \underline{Proof}: Note that 
\begin{equation}\label{energy difference}
 E(f_{min},s) = E(f_{min},0) + \int\limits_{\mathbb{R}} \left\{ \frac{1}{2}(s')^2 + \frac{\lambda_\sigma}{4}(s^2 - 2)s^2 + \frac{\beta}{2} f^2 s^2 + \frac{d^2}{2 R^2} s^2\right\}
\end{equation}
If we can show that the second term is negative we'd be done. To this end, take $$s = \frac{1}{\cosh(Bx)}$$ where $B = \sqrt{\frac{\lambda_\phi}{2}}$. Thus, plugging $s$ into the second term we have that
\begin{eqnarray*}
 & & \int\limits_{\mathbb{R}} \left\{ \frac{1}{2}(s')^2 + \frac{\lambda_\sigma}{4}(s^2 - 2)s^2 + \frac{\beta}{2} f^2 s^2 + \frac{d^2}{2 R^2} s^2\right\} \\
 &=& \int\limits_{\mathbb{R}}\left[  \frac{B^2}{2} - \frac{\lambda_\sigma}{4} + \frac{\beta}{2}   \right] \frac{\sinh^2(B y^1)}{\cosh^4(B y^1)} + \int\limits_{\mathbb{R}} \left[\frac{d^2 }{2 R^2} - \frac{\lambda_\sigma}{4} \right] \frac{1}{\cosh^2(B y^1)} \\
 &=& \frac{1}{B}\left\{ \frac{\beta}{3} + \frac{d^2}{R^2} - \frac{\lambda}{2} \right\}
\end{eqnarray*}
The additional constraint (\ref{non-zero s minimizer further constraint}) implies that $E(f_{min},s) < E(f_{min},0)$ for sufficiently large $R$.
\begin{flushright}
 $\Box$
\end{flushright}

\subsubsection{Interface Evolution and Current Quenching}\label{current quenching}
In this section we make two observations about the approximate solution $F_0$ and the surface $R$ about which it is concentrated.
\bigskip

The first observation we make is that when the interface we find has a current (i.e. $s_0\neq 0$), then the interface moves towards the origin. To see why this is true, recall that the surface $R$ satisfies the geometric relation
\[
 \frac{1}{\sqrt{1 - (R')^2}} \left( \frac{R''}{1 - (R')^2} + \frac{1}{R}\right) = \frac{1}{\sqrt{1 - (R')^2}} \frac{d^2}{R^3} \frac{\left\|s_0\right\|_{L^2(\mathbb{R})}^2}{\left\| F_0'\right\|_{L^2(\mathbb{R})}^2}
\]
Rearranging, this can be restated as
\[
 R'' = \left[\frac{d^2}{R^2} \frac{\left\|s_0\right\|_{L^2(\mathbb{R})}^2}{\left\| F_0'\right\|_{L^2(\mathbb{R})}^2} - 1\right]\frac{1 - (R')^2}{R}
\]
Since $F_0$ minimizes $\mu(f,s;R)$, then $F_0$ satisfies
\[
 -F_0'' + \nabla_\Phi W(F_0,R) = 0
\]
Multiplying this by $F_0'$ and integrating, we also have that $F_0$ satisfies
\[
 \frac{1}{2}(F_0')^2 = W(F_0,R)
\]
Thus,
\begin{equation}\label{R rearranged}
 R'' = \left[\frac{1}{2} \frac{d^2}{R^2} \frac{\left\|s_0\right\|_{L^2(\mathbb{R})}^2}{\left\| W(F_0,R)\right\|_1} - 1\right]\frac{1 - (R')^2}{R}
\end{equation}
Since $$W(F_0,R) = V(F_0) + \frac{d^2}{2 R^2} s_0^2 > \frac{d^2}{2 R^2} s_0^2$$ as $V(F_0) \geq 0$, then $$\frac{1}{2} \frac{d^2}{R^2} \frac{\left\|s_0\right\|_{L^2(\mathbb{R})}^2}{\left\| W(F_0,R)\right\|_1} > 1$$ and so $R'' < 0$ (as long as $\left|R'\right| < 1$). Since $R'(0) = 0$, this implies that $R$ is moving towards the origin.
\begin{figure}[H]
 \centering
 \includegraphics[scale=0.4]{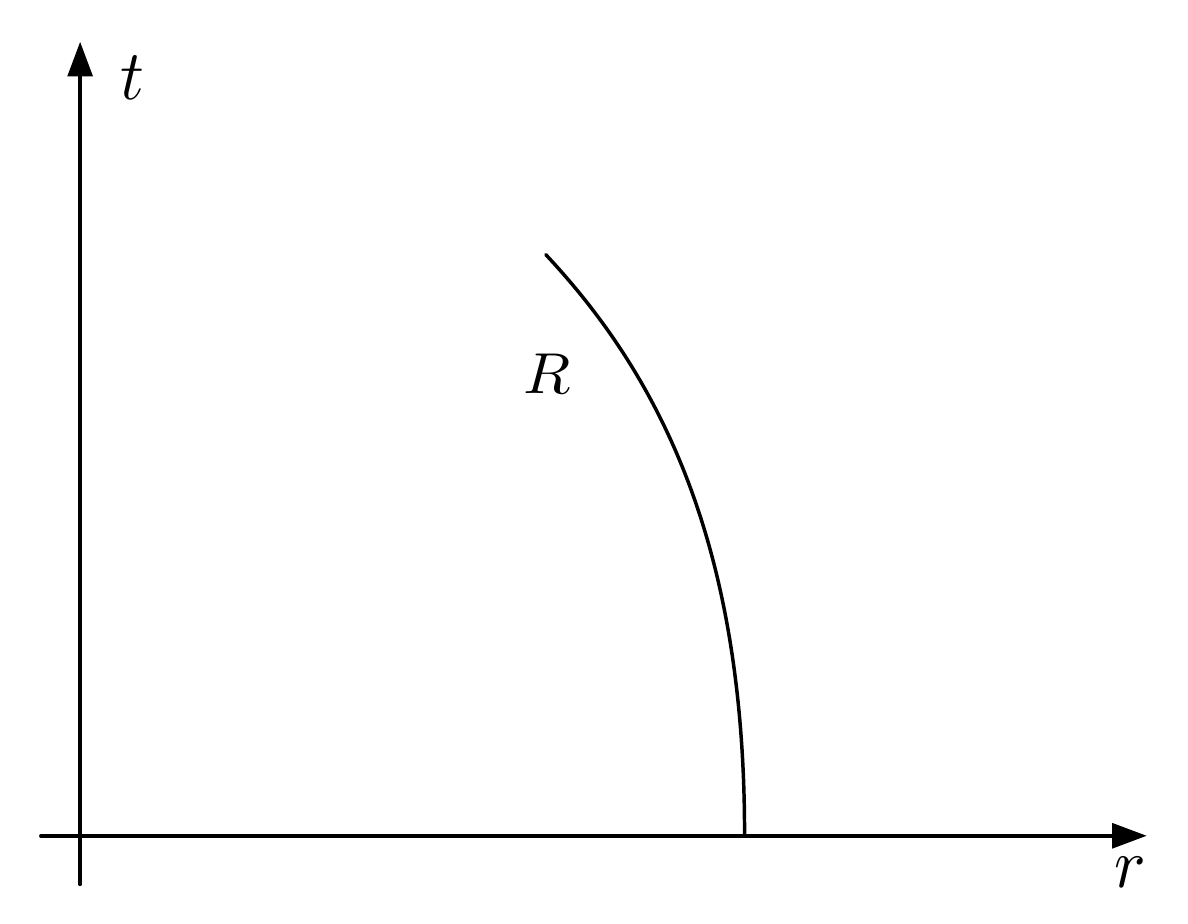}\label{R decreasing}
 \caption{Since $R'(0) = 0$ and $R''< 0$ whenever $s_0\neq 0$ then $R$ is decreasing, at least for a short time.}
\end{figure}
\bigskip

The second observation we make is that for $R$ sufficiently small, then for $F = (f,s)$ minimizing $\mu(f,s;R)$ we necessarily have that $s=0$. By lemma \ref{minimizer properties}, for $(f,s)$ minimizing $\mu(f,s;R)$, then $\left| f\right| \leq 1$ and $\left|s\right| \leq 1$. Then,
\[
 W(f,s;R) = W(f,0;R) + \left[\int\limits_0^1 \partial_\sigma V(f,\lambda s) d\lambda + \frac{d^2}{R^2}\right] s^2
\]
By condition 1 of (\ref{potential assumptions}), $\left|\partial_\sigma V(f,\lambda s)\right| \leq c$ for all $0\leq \lambda \leq 1$. Thus, for sufficiently small $R > 0$ $$W(f,s;R) > W(f,0;R)$$ This implies that for sufficiently small $R$, then $(f,s)$ minimizing $\mu(f,s;R)$ implies that $s=0$.
\bigskip

Suppose $(f,s)$ minimize $\mu(f,s;R)$ and suppose we have a potential $V$ for which there exists a range of $R$ so that $s(y^1;R)\neq 0$. By the second observation, we see that even though $s(y^1;R)\neq 0$ for some $R$, there exists $R$ sufficiently small for which $s(y^1;R) = 0$. Suppose $R_*$ is the smallest value of $R$ so that $s(y^1;R_* + \delta)\neq 0$ for  $0< \delta \ll 1$. By the first observation, if $R_* < R(0) < R_* + \delta$, then $R$ becomes smaller as the system evolves. Our solution only makes sense when $\left| R'\right| < 1$. So the question becomes, does $R$ become smaller than $R_*$ before $\left| R'\right| =1$ and hence do we have current quenching?
\bigskip

Pick $R_* < R(0) = R_* + \delta$. As we observed before, $R$ becomes smaller under the flow of (\ref{R}). Pick $\delta$ small enough so that 
\[
 -1 \leq \frac{1}{2} \frac{d^2}{R^2} \frac{\left\|s_0\right\|_{L^2(\mathbb{R})}^2}{\left\| W(F_0,R)\right\|_1} - 1 \leq -\frac{1}{2}
\]
for $\frac{1}{2} R_* \leq R \leq R_* + \delta$ (for $R<R_*$ this quantity is exactly $-1$ as $s=0$). Further, when $\frac{1}{2} R_* \leq R \leq R_* + \delta$, then we can estimate (\ref{R rearranged}) as
\[
 -2 \frac{1 - (R')^2}{R_*} \leq R'' \leq -\frac{1}{2}\frac{1 - (R')^2}{R_* + \delta}
\]
It follows (after integrating) that
\begin{equation}\label{R' estimate}
 -\tanh(\frac{4}{R_*} y^0) \leq R' \leq -\tanh(\frac{4}{R_*} y^0)
\end{equation}
From this, we can see that $\left| R'\right| < 1$ as long as $\frac{1}{2} R_* \leq R \leq R_* + \delta$. Furthermore, we can integrate once more in $y^0$ to find that
\begin{equation}\label{R estimate}
 (R_* + \delta) - \frac{R_*}{4} \log\cosh(\frac{4}{R_*}y^0) \leq R \leq (R_* + \delta)\left[1 - \log\cosh(\frac{y^0}{R_* + \delta})\right]
\end{equation}
again as long as $\frac{1}{2} R_* \leq R \leq R_* + \delta$. From (\ref{R estimate}) and (\ref{R' estimate}) we see that for our system $R$ becomes smaller than $R_*$ in finite time. Thus, solutions we find undergo current quenching.

\subsection[Asymptotics of F0 and F1]{Asymptotics of $F_0$ and $F_1$}
 We would like to examine $F_0$ for large values of $y^1$. This leads us to the following proposition
\begin{proposition}\label{asymptotics}
 Suppose $(f,s)\in\mathcal{A}$ solves the minimization problem (\ref{minimization problem}). Then there exists $\alpha > 0$ so that
 \begin{equation}\label{f and s minimizer asymptotic}
  \begin{array}{ccc} \left\{\begin{array}{cl} 1 - \left| f\right| &\lesssim e^{-\alpha \left|y^1\right|} \\ \left| f'\right| &\lesssim e^{-\alpha \left|y^1\right|} \end{array}\right. & \textrm{ and } & \left\{\begin{array}{cl} \left| s\right| &\lesssim e^{-\alpha \left| y^1\right|} \\ \left|s'\right| &\lesssim e^{-\alpha \left| y^1\right|} \end{array}\right. \end{array}
 \end{equation}
\end{proposition}
Minimizers $F = (f,s)$ of (\ref{minimization problem}) satisfy $-F'' + \nabla_\Phi W(F,R) = 0$. Using assumption 3 of (\ref{potential assumptions}) one can easily obtain (\ref{f and s minimizer asymptotic}).
\bigskip

\begin{proposition}\label{asymptotics of F_1}
 Suppose $F_1$ solves (\ref{F_1 equation}). Then there exists $\alpha > 0$ so that for $\beta = 0,1$ we have
 \begin{equation}\label{F_1 asymptotics}
 \begin{array}{ccc}
  \left| \partial_{y^1}^\beta F_1\right| &\lesssim& e^{-\alpha \left| y^1\right|}
  \end{array}
 \end{equation}
\end{proposition}
Since the left hand side of (\ref{F_1 equation}) decays exponentially fast, again using assumption 3 of (\ref{potential assumptions}) and using standard arguments yields (\ref{F_1 asymptotics}).

\section{Effective Dynamics}
\subsection{Approximation Using Profiles Coming From the Formal Asymptotics}
The main question we would like to answer is: Suppose $\Phi$ is a solution to (\ref{Wave equation in minkowski-normal coordinates}) with the following properties
\begin{itemize}
 \item $\Phi = (-1,0) \textrm{ for } y_1 < -y_*^1 \;\; \text{ and } \;\; \Phi = (1,0) \textrm{ for } y^1 > y_*^1$
 
 \item $\Phi(0,y^1) \textrm{ is close to } F_0(\frac{y^1}{\epsilon};R(0)) + \epsilon F_1(\frac{y^1}{\epsilon};R(0),R'(0))$
 
 \item $\partial_{y^0}\Phi(0,y^1) \textrm{ is close to } \partial_{y^0} \bigg|_{y^0 = 0} ( F_0(\frac{y^1}{\epsilon};R) + \epsilon F_1(\frac{y^1}{\epsilon};R,R'))$
\end{itemize}
then does $\Phi$ remain close to $F_0(\frac{y^1}{\epsilon};R) + \epsilon F_1(\frac{y^1}{\epsilon};R,R')$ up to some $y_*^0$ independent of $\epsilon$?
\bigskip

We will use the translation symmetry of the profiles $F_0$ to find a function $a:[0,y_*^0]\rightarrow\mathbb{R}$ so that for each $y^0\in [0,y_*^0]$ the difference between a solution $\Phi$ to (\ref{superconducting interface model in minkowski normal coordinates}) and $F_0(\frac{y^1 - a(y^0)}{\epsilon};R(y^0))$ is minimized. For each $R\in\mathbb{R}$ and $\Psi\in L^2(\mathbb{R},\mathbb{R}^2)$, define
\begin{equation}\label{h equation}
 h_{\epsilon}(\Psi,R,a) := \left\| \Psi - F_0(\frac{y^1 - a}{\epsilon};R)\right\|_{L^2(\mathbb{R})}^2
\end{equation}
\begin{equation}\label{first derivative of h}
 G_{\epsilon}(\Psi,R,a) := \partial_a h_{\epsilon}(\Psi,R,a) = -\frac{2}{\epsilon} \left<\Psi - F_0(\frac{y^1 - a}{\epsilon};R), \partial_{y^1} F_0(\frac{y^1 - a}{\epsilon};R) \right>_{L^2(\mathbb{R})}
\end{equation}
For each $y^0\in [0,y_*^0]$ we want to find a sufficiently regular $a(y^0)$ so that $$G_{\epsilon}(\Phi(y^0),R(y^0),a(y^0)) = 0$$

Define
\begin{equation}\label{U delta neighbourhood}
 U_{\delta,\epsilon} := \left\{ (\Psi,R)\in L^2(\mathbb{R},\mathbb{R}^2)\times \mathbb{R} \; : \; \inf\limits_{a\in\mathbb{R}} h_{\epsilon}(\Psi,R,a) < \delta\right\}
\end{equation}
\begin{equation}\label{V delta neighbourhood}
 V_{\delta,\epsilon}(a_0) := \left\{ (\Psi,R) \in L^2(\mathbb{R},\mathbb{R}^2)\times\mathbb{R} \; : \; h_{\epsilon}(\Psi,R,a_0) < \delta \right\}
\end{equation}
\begin{lemma}\label{existence of tilde a}
 There exists $\delta > 0$ and a unique $C^1$ map, $C^1$ with respect to the $L^2\times\mathbb{R}$ topology, $\tilde{a}: U_{\delta,\epsilon}\rightarrow\mathbb{R}$ so that $G_{\epsilon}(\Psi,R,\tilde{a}(\Psi,R))=0$ where both $\delta$ and $\tilde{a}$ possibly depend on $\epsilon$.
\end{lemma}

\noindent\underline{Proof}: Since
\begin{enumerate}
 \item $G_{R,\epsilon}$ is $C^1$ in $\Psi$ as it is linear in $\Psi$, is $C^1$ in $R$ as $F_0$ is $C^1$ in $R$, and is $C^1$ in $a$ because $F_0$ is $C^2$ in $y^1$
 
 \item $G_{\epsilon}(F_0(\frac{\cdot - a_0}{\epsilon};R),R,a_0) = 0$
 
 \item $$\partial_a \big|_{a = a_0} G_{\epsilon}(F_0(\frac{\cdot - a_0}{\epsilon};R),R,a) = \left<\partial_{y^1} F_0(\frac{\cdot - a_0}{\epsilon};R),\partial_{y^1} F_0(\frac{\cdot - a_0}{\epsilon};R)\right>_{L^2(\mathbb{R})} > 0$$
\end{enumerate}
then we can apply the implicit function. That is, there exists $\delta > 0$ and a unique $C^1$ map $a:V_{\delta,\epsilon}(a_0)\rightarrow\mathbb{R}$, both $\delta$ and $a$ possibly depending on $\epsilon$, so that $G_{\epsilon}(\Psi,R,a(\Psi,R)) = 0$ for all $(\Psi,R)\in V_{\delta,\epsilon}(a_0)$.
\bigskip

Observe that $$ U_{\delta,\epsilon} = \bigcup_{b\in\mathbb{R}} V_{\delta,\epsilon}(a_0+b)$$ For each $(\Psi,R)\in U_{\delta,\epsilon}$ there exists $b\in\mathbb{R}$ so that for $\tau_b \Psi := \Psi(\cdot - b)$, then $(\tau_b\Psi,R)\in V_{\delta,\epsilon}(a_0)$. Define $\tilde{a}_b(\Psi,R) = a(\tau_b \Psi,R) + b$, then $G_{\epsilon}(\Psi,R,\tilde{a}_b(\Psi,R)) = 0$. If $\tau_c \Psi$, $\tau_b\Psi\in V_{\delta,\epsilon}(a_0)$, then by the uniqueness of $a$ one has that $\tilde{a}_b(\Psi,R) = \tilde{a}_c(\Psi,R) + c - b$ and thus $a(\tau_b \Psi,R) = a(\tau_c \Psi,R)$. Therefore, one can find a unique $\tilde{a}: U_{\delta,\epsilon}\rightarrow\mathbb{R}$ so that $G_{\epsilon}(\Psi,R,\tilde{a}(\Psi,R)) = 0$.
\begin{flushright}
 $\Box$
\end{flushright}

Suppose $\Phi$ is a solution to (\ref{equations of motion in polar coordinates}) and that $(\Phi(0),R(0))\in U_{\delta,\epsilon}$ with $\delta > 0$ coming from lemma \ref{existence of tilde a}. Then, there exists some maximal $0 < y_*^0(\epsilon)$, where $y_*^0(\epsilon)$ may or may not depend on $\epsilon$ as $U_{\delta,\epsilon}$ depends on $\epsilon$, so that $(\Phi(y^0),R(y^0))\in U_{\delta,\epsilon}$ for all $0 \leq y^0 \leq y_*^0(\epsilon)$. Thus, $$G_{\epsilon}(\Phi(y^0),R(y^0),\tilde{a}(\Phi(y^0),R(y^0))) = 0$$ for $0 \leq y^0 \leq y_*^0(\epsilon)$. While proving theorem \ref{main theorem intro}, we will actually show that $y_*^0(\epsilon)$ does not depend on $\epsilon$. For if $y_*^0(\epsilon)$ does depend on $\epsilon$, then $(\Phi(y_*^0(\epsilon),R(y_*^0(\epsilon))\in U_{\delta,\epsilon}$ which contradicts the maximality of $y_*^0(\epsilon)$.
\begin{corollary}\label{existence of a}
 Suppose $\Phi$ solves (\ref{equations of motion in polar coordinates}) and suppose that $(\Phi(0),R(0))\in U_{\delta,\epsilon}$ with $\delta > 0$ from lemma \ref{existence of tilde a}. Then, there exists $y_*^0(\epsilon)>0$ and a unique $C^1$ function $a(y^0)$ so that $G_\epsilon(\Phi(y^0),R(y^0),a(y^0)) = 0$ for all $0\leq y^0 \leq y_*^0(\epsilon)$.
\end{corollary}

Given a solution $\Phi$ to (\ref{equations of motion in polar coordinates}), we define 
\begin{equation}\label{F_0 tilde}
 \tilde{F}_0(y^0,y^1) = F_0(\frac{y^1 - a(y^0)}{\epsilon};R(y^0))
\end{equation}
where $a(y^0)$ is from corollary \ref{existence of a}.
\bigskip

Examining (\ref{F_1 equation}) we see that $F_1$ does \underline{not} have a translation symmetry in $y^1$ as the inhomogeneity of (\ref{F_1 equation}) depends explicitly on $y^1$. Instead, we have that for $v = R'$, then $F_1(\frac{y^1 - a}{\epsilon};R,v)$ from proposition \ref{existence of F_1} solves
\begin{equation}\label{translated F_1 equation}
 \epsilon L_\epsilon(\tilde{F}_0;R)\left[ F_1(\frac{y^1 - a}{\epsilon};R,R')\right] = H(R) \partial_{y^1}\tilde{F}_0 - \frac{1}{\epsilon} m(y^0) (\frac{y^1 - a}{\epsilon}) \partial_R w(\tilde{F}_0,R) 
\end{equation}
where $m$ was defined (\ref{m and n}). Remember that $F_1$ was defined independent of $\epsilon$. We define
\begin{equation}\label{F_1 tilde}
 \tilde{F}_1(y^0,y^1) = F_1(\frac{y^1 - a(y^0)}{\epsilon};R(y^0),R'(y^0))
\end{equation}

For our result, we'll need to control two quantities. The first is the error and the second is the \textbf{shift} $a(y^0)$. Define the \textbf{error} between $\Phi$ and our approximation $\tilde{F}_0 + \epsilon \tilde{F}_1$ as
\begin{equation}\label{error}
 \xi := \Phi - \tilde{F}_0 - \epsilon \tilde{F}_1
\end{equation}
and define the quantity
\begin{equation}\label{shift}
 \underline{A}(y^0) := \left( 1 + \frac{\left| a(y^0) \right|}{\epsilon} + \frac{\left| a'(y^0) \right|}{\epsilon} \right)^3
\end{equation}
\bigskip

An observation that we will make use of later is the following. Since $0 = \partial_a h(a(y^0))$ and $\tilde{F}_1\perp \partial_{y^1} F_0$, we can use (\ref{first derivative of h}) to get that 
\begin{equation}\label{xi orth F0'}
 0 = \int\limits_\mathbb{R} \xi\cdot \partial_{y^1} \tilde{F}_0
\end{equation}
where we needed to use the fact that $\tilde{F}_1\perp \partial_{y^1} \tilde{F}_0$, from proposition \ref{existence of F_1}, to go from the second line to the third. That is, we have that $\xi\perp \partial_{y^1} \tilde{F}_0$.
\bigskip

Next, we will plug $\Phi = \tilde{F}_0 + \epsilon \tilde{F}_1 + \xi$ into (\ref{Wave equation in minkowski-normal coordinates}) and find the equation that $\xi$ solves. Doing so, we find that $\xi$ solves
\begin{equation}\label{xi equation}
\frac{m^2}{n^2} \partial_{y^0 y^0} \xi + B^\alpha \partial_\alpha\xi + L_\epsilon(\tilde{F}_0,R)\xi + S_{-1} + S_0 + N = 0
\end{equation}
where we used the fact that $$-\partial_{y^1 y^1} \tilde{F}_0 + \frac{1}{\epsilon^2} w(\tilde{F}_0,R) = 0$$ to simplify and we defined
\begin{equation}\label{S_-1}
 S_{-1} = \epsilon L_\epsilon(\tilde{F}_0,R) \tilde{F}_1 + B^1 \partial_{y^1} \tilde{F}_0 + \frac{y^1}{\epsilon^2} m(y^0) \partial_R w(\tilde{F}_0,R) 
\end{equation}
\begin{equation}\label{S_0}
 S_0 = \frac{m^2}{n^2} \partial_{y^0 y^0}(\tilde{F}_0 + \epsilon \tilde{F}_1) + B^0 \partial_{y^0} (\tilde{F}_0 + \epsilon \tilde{F}_1) + \epsilon B^1 \partial_{y^1} \tilde{F}_1
\end{equation}
\begin{equation}\label{N}
 N = \frac{1}{\epsilon^2}\left[ w(\tilde{F}_0 + \tilde{F}_\xi,R + y^1 m) - w(\tilde{F}_0,R) - \operatorname{Hess}_\Phi W(\tilde{F}_0,R)\tilde{F}_\xi - y^1 m(y^0) \partial_R w(\tilde{F}_0,R)  \right]
\end{equation}
where we've defined $\tilde{F}_\xi = \epsilon \tilde{F}_1 + \xi$. Note that (\ref{xi equation}) only makes sense on $(0,y_*^0)\times (-y_*^1,y_*^1)$, but since $\Phi$ and $\tilde{F}_0 + \epsilon \tilde{F}_\xi$ are both defined on $(0,y_*^0)\times \mathbb{R}$, then $\xi$ is defined on this set too. Recall that we have that $\Phi = (1,0)$ for $y^1 > y_*^1$ and $\Phi = (-1,0)$ for $y^1 < -y_*^1$. Outside of $\left| y^1 \right| \leq y_*^1$, we use the asymptotics derived in proposition \ref{asymptotics} to get that
\begin{gather}\label{xi estimate outside}
 \left\| \xi\right\|_{L^2({\left| y^1 \right| > y_*^1})} 
 \leq \left\| (\operatorname{sgn}(y^1),0) - \tilde{F}_0 \right\|_{L^2({\left| y^1 \right| > y_*^1})} +  \epsilon \left\| \tilde{F}_1 \right\|_{L^2({\left| y^1 \right| > y_*^1})}
 \lesssim e^{- \alpha \frac{y_*^1 - a}{\epsilon}} 
\end{gather}
 for some $\alpha > 0$. Thus, we have that $\xi$ is small outside of $\left| y^1 \right| \leq y_*^1$ if we can control the size of $a$ and if $\epsilon$ is taken sufficiently small. We are then left to estimate $\xi$ on $\left| y^1 \right| \leq y_*^1$. We use the following quantities to control $\xi$ on $\left| y^1 \right| \leq y_*^1$.
\begin{definition}
  For $Q = (Q_1,Q_2): (0,y_*^0)\times (-y_*^1,y_*^1) \rightarrow \mathbb{R}^2$ define the energy density
 \begin{equation}\label{energy density}
  e(Q) = \frac{1}{2} \frac{m^2}{n^2} \partial_{y^0} Q^2 + \frac{1}{2}\partial_{y^1} Q^2 + \frac{1}{2\epsilon^2} Q\cdot \operatorname{Hess}_\Phi W(\tilde{F}_0,R)Q
 \end{equation}
  Using the energy density, we define the energy of $Q$ as
 \begin{equation}\label{energy}
  E(Q) = \int\limits_{\left| y^1 \right| \leq y_*^1} e(Q) dy^1  
 \end{equation}
For convenience we set $E(y^0) = E(\xi)(y^0)$.
\end{definition}

Using this new definition, we obtain a very useful corollary to theorem \ref{spectral estimate theorem} that we will use to control the error term $\xi$
\begin{corollary}\label{energy-spectral estimate corollary}
 Suppose $\Phi$ is a solution to (\ref{Wave equation in minkowski-normal coordinates}) with initial data as described in section \ref{Initial Data and the Existence of Solutions}. For $\xi$ as defined in (\ref{error}), we have that
 \begin{equation}\label{energy-spectral estimate}
 \frac{1}{\epsilon^2} \int\limits_{-y_*^1}^{y_*^1} \left| \xi \right|^2 \lesssim E + \frac{1}{\epsilon^2} e^{-\alpha \frac{y_*^1 - a}{\epsilon}}
\end{equation}
for some $\alpha > 0$
\end{corollary}
\noindent \underline{Proof of corollary \ref{energy-spectral estimate corollary}}: Use theorem \ref{spectral estimate theorem} along with (\ref{xi estimate outside}) to obtain the estimate (\ref{energy-spectral estimate}).
\begin{flushright}
 $\Box$
\end{flushright}

\subsection{Main Result}
The main theorem of this paper is the following.
\begin{theorem}\label{main theorem}
 Suppose $\Phi$ solves (\ref{Wave equation in minkowski-normal coordinates}). Further assume that $\Phi(0)$, $\partial_{y^0} \Phi(0)$, $a(0)$, and $a'(0)$ satisfy
 \begin{equation*}
  \begin{array}{ccc} \underline{A}(0) \lesssim 1 & \textrm{ and } & E(0) \lesssim \epsilon^2 \end{array}
 \end{equation*}
 Then there exists $0 < \overline{y}^0 < y_*^0$, $\overline{y}^0$ independent of $\epsilon$, and $a:(0,\overline{y}^0)\rightarrow \mathbb{R}$ so that
 \begin{equation*}
  \begin{array}{ccc} \underline{A}(y^0) \lesssim 1 & \textrm{ and } & E(y^0) \lesssim \epsilon^2 \end{array}
 \end{equation*}
 for all $0 \leq y^0 \leq \overline{y}^0$.
\end{theorem}

To prove theorem \ref{main theorem}, we will use the following two estimates
\begin{theorem}[Energy Estimate Theorem]\label{energy estimate theorem}
 Suppose $\Phi$ solves (\ref{Wave equation in minkowski-normal coordinates}). Then for as long as $a(y^0)$ is well defined we have
 \begin{align}\label{energy estimate}
  & \left[1 - \frac{y_*^1}{(R - y_*^1 m)^3}\right]E - E(0) \nonumber \\
  \lesssim &  \; \left[ \sqrt{\epsilon}^{3} \sqrt{E} \underline{A} + \sqrt{\epsilon} \sqrt{E}^{3} + \sqrt{\epsilon} e^{-\alpha \frac{y_*^1 - a}{\epsilon}} \underline{A} + \frac{1}{\epsilon} \sqrt{E} e^{-\alpha \frac{y_*^1 - a}{\epsilon}} \right]_0^{y^0} + \int\limits_0^{y^0} (E + \epsilon \sqrt{E} + \frac{1}{\epsilon^2} e^{-\alpha \frac{y_*^1 - a}{\epsilon}})\underline{A}
 \end{align}
\end{theorem}
\begin{theorem}[Bounded Shift Theorem]\label{bounded shift theorem}
 Suppose $\Phi$ solves (\ref{Wave equation in minkowski-normal coordinates}). Then for as long as $a(y^0)$ is well defined we have
 \begin{align}\label{bounded shift}
  \left[1 -  \epsilon \underline{A} - \sqrt{\epsilon} \underline{A} \sqrt{E} - e^{-\alpha \frac{y_*^1 - a}{\epsilon}} \right] \frac{\left| a'' \right|}{\epsilon}  \; \lesssim &  \; \frac{1}{\sqrt{\epsilon}^{5}} (\epsilon\underline{A} + \sqrt{(\epsilon \sqrt{E} + e^{-\alpha \frac{y_*^1 - a}{\epsilon}})} \sqrt[4]{E})(\sqrt{\epsilon}^{3} + \epsilon \sqrt{E} + e^{-\alpha \frac{y_*^1 - a}{\epsilon}})
 \end{align}
\end{theorem}
\bigskip

Theorem \ref{main theorem intro} is obtained from theorem \ref{main theorem} by applying the spectral estimate (\ref{energy-spectral estimate}) and using the estimate $E(y^0) \lesssim \epsilon^2$. Assuming theorems \ref{energy estimate theorem} and \ref{bounded shift theorem} are true, we can prove theorem \ref{main theorem}.
\bigskip

\noindent \underline{Proof of Theorem \ref{main theorem}}: We will be implementing a bootstrap argument to prove theorem \ref{main theorem}. In order to close the argument to be outlined, we may need to choose $y_*^0$ and $y_*^1$ smaller, still independent of $\epsilon$, so that
\begin{gather}\label{y_*^1 suff small}
 C \frac{y_*^1}{(R - y_*^1 m)^3} \leq \frac{1}{2}
\end{gather}
where $C = C(y_*^0, y_*^1)$ is the constant from theorems \ref{energy estimate theorem} and \ref{bounded shift theorem}. We can find such a $y_*^1$, because $C(y_*^0,y_*^1)\rightarrow 0$ as $y_*^1\rightarrow 0$. We'll make use of the following two estimates in order to complete the proof
\begin{gather}\label{a estimates}
 \left| a\right| \leq \left|a(0)\right| + \int\limits_0^{y^0} \left| a'\right| \;\;\;\;\;  \text{ and } \;\;\;\;\;  \left| a'\right| \leq \left|a'(0)\right| + \int\limits_0^{y^0} \left| a''\right|
\end{gather}

Next, suppose $a(y^0)$ is well defined on the interval $I = (0,b)$. Define 
\[
 \begin{array}{ccc} E_M(I) := \max\limits_{y^0\in I} E(y^0) & \textrm{ and } & \underline{A}_M(I) := \max\limits_{y^0\in I} \underline{A}(y^0)  \end{array}
\]
Using theorem \ref{energy estimate theorem} we have that
\begin{align}\label{max energy estimate}
 E_M(I) \; & \lesssim \; E(0) + \frac{1}{2} E_M(I) + \epsilon^{3/2} \underline{A}_M(I) E_M(I)^{1/2} + \epsilon^{1/2} E_M(I)^{3/2} + \epsilon^{1/2} e^{-\frac{\alpha}{\epsilon} y_*^1} e^{\alpha \underline{A}_M(I)} \\
 & \; + \frac{1}{\epsilon} E_M(I)^{1/2} e^{-\frac{\alpha}{\epsilon} y_*^1} e^{\alpha \underline{A}_M(I)} + \left| I\right| \left[ E_M(I) + \epsilon E_M(I)^{1/2} + \frac{1}{\epsilon^2} e^{-\frac{\alpha}{\epsilon} y_*^1} e^{\alpha \underline{A}_M(I)}\right] A_M(I) \nonumber
\end{align}
Using theorem \ref{bounded shift theorem} and (\ref{a estimates}) we have that
\begin{align}\label{max bounded shift estimate}
 \underline{A}_M(I)^{1/3} \;& \lesssim \; (1 + \left| I\right|) \underline{A}(0)^{1/3} + \frac{\left| I\right| + \left| I\right|^2}{1 - \epsilon \underline{A}_M(I) - \sqrt{\epsilon} E_M(I)^{1/2} \underline{A}_M(I) + e^{-\frac{\alpha}{\epsilon} y_*^1 } e^{\alpha \underline{A}_M(I)} } \underline{B}_M(I) 
\end{align}
where we've introduced $\underline{B}_M(I)$ as
\begin{align}\label{intermediate quantity}
 \underline{B}_M(I) &:=& \frac{1}{\sqrt{\epsilon}^{5}} \left[ \epsilon  \underline{A}_M(I) +\sqrt{\epsilon \sqrt{E_M(I)} + e^{-\frac{\alpha}{\epsilon} y_*^1} e^{\alpha \underline{A}_M(I)}} \sqrt[4]{E_M(I)}\right]\left[\sqrt{\epsilon}^{3} + \epsilon \sqrt{E_M(I)} + e^{-\frac{\alpha}{\epsilon} y_*^1} e^{\alpha\underline{A}_M(I)} \right]
\end{align}
By corollary \ref{existence of a}, $a(y^0)$ is well defined up to some time $\widehat{y}^0$ (corollary \ref{existence of a} doesn't tell us that $\widehat{y}^0$ is independent of $\epsilon$, just that it exists). There exists $\overline{y}^0$ so that for $I = \left(0, \min\left\{ \widehat{y}^0, \overline{y}^0 \right\}\right)$, then estimates (\ref{max energy estimate}) and (\ref{max bounded shift estimate}) imply that
\[
 \begin{array}{ccc} E_M(I) \lesssim \epsilon^2 & \textrm{ and } & \underline{A}_M(I) \lesssim 1 \end{array}
\]
\smallskip

If $\min\left\{ \widehat{y}^0,\overline{y}^0\right\} = \overline{y}^0$, then we are done. If not, then because $$\underline{E}_M(I) \lesssim \epsilon^2$$ using corollary \ref{existence of a} we actually have that $a(y^0)$ exists beyond $\widehat{y}^0$. Boot strapping allows us to conclude that $a(y^0)$ exists and is well defined on $I = \left(0,\overline{y}^0\right)$ and on $I$ that
\[
 \begin{array}{ccc} E_M(I) \lesssim \epsilon^2 & \textrm{ and } & \underline{A}_M(I) \lesssim 1 \end{array}
\]
\begin{flushright}
 $\Box$
\end{flushright}

\subsection{Proof of Energy Estimate (Theorem \ref{energy estimate theorem})}
We require an estimate of $(y^1)^\gamma \partial_{y^1}^\alpha \partial_{y^0}^\beta \tilde{F}_i$, for $\alpha$, $\beta$, $\gamma \in \mathbb{N}\cup \left\{0\right\}$, to prove this theorem and theorem \ref{bounded shift theorem}. A point on notation before continuing. We have that $F_1 = F_1(y^1;R,v)$ where third slot of $F_1$ is called $v$.
\begin{lemma}\label{partial F_i estimate lemma}
 For $\alpha$, $\beta$, $\kappa$, $\lambda$, and $\gamma \in \mathbb{N}\cup \left\{0\right\}$, then for $(R,v)\in (r_0,r_1)\times (-1,1)$, $r_0 < r_1$ coming from the non-degeneracy condition (\ref{nondeg con}), and $a\in \left\{ a(y^0) \; : \; 0\leq y^0 \leq y_*^0  \right\}$ we have that 
 \begin{equation}\label{partial F_0 estimate}
  \left\|(y^1)^\gamma (\frac{\partial}{\partial y^1})^\alpha (\frac{\partial}{\partial R})^\beta \left[F_0(\frac{\cdot - a}{\epsilon};R)\right] \right\|_{L^2(\mathbb{R})} \lesssim \frac{1}{\epsilon^{\alpha - \gamma - \frac{1}{2}}} \left( 1 + \frac{\left| a(y^0)\right|}{\epsilon} \right)^\gamma
 \end{equation}
 \begin{equation}\label{partial F_1 estimate}
  \left\|(y^1)^\gamma (\frac{\partial}{\partial y^1})^\alpha (\frac{\partial}{\partial R})^\beta (\frac{\partial}{\partial v })^\kappa \left[F_1(\frac{\cdot - a}{\epsilon};R,v)\right] \right\|_{L^2(\mathbb{R})} \lesssim \frac{1}{\epsilon^{\alpha - \gamma - \frac{1}{2}}} \left( 1 + \frac{\left| a(y^0)\right|}{\epsilon} \right)^\gamma
 \end{equation}

 where the constant in the estimate depends on $y_*^0$, but not $\epsilon$.
\end{lemma}
\noindent

\noindent \underline{Proof of lemma \ref{partial F_i estimate lemma}}: For (\ref{partial F_0 estimate}) we have
\begin{eqnarray*}
 \left\|(y^1)^\gamma (\frac{\partial}{\partial y^1})^\alpha (\frac{\partial}{\partial R})^\beta \left[F_0(\frac{\cdot - a}{\epsilon};R)\right] \right\|_{L^2(\mathbb{R})} 
 &=& \frac{1}{\epsilon^\alpha} \left\|(y^1)^\gamma \frac{\partial^{\alpha + \beta} F_0}{\partial (y^1)^\alpha \partial R^\beta}(\frac{\cdot - a}{\epsilon};R) \right\|_{L^2(\mathbb{R})} \\
 &\lesssim& \frac{1}{\epsilon^{\alpha - \gamma}} \left\|\left(\frac{y^1 - a + a}{\epsilon}\right)^\gamma \frac{\partial^{\alpha + \beta} F_0}{\partial (y^1)^\alpha \partial R^\beta}(\frac{\cdot - a}{\epsilon};R) \right\|_{L^2(\mathbb{R})} \\
 &\lesssim& \frac{1}{\epsilon^{\alpha - \gamma}} \left\| \left[ \sum\limits_{j=0}^\gamma {\gamma \choose j} \left(\frac{y^1 - a}{\epsilon}\right)^{\gamma - j} \left(\frac{a}{\epsilon}\right)^j \right] \frac{\partial^{\alpha + \beta} F_0}{\partial (y^1)^\alpha \partial R^\beta}(\frac{\cdot - a}{\epsilon};R) \right\|_{L^2(\mathbb{R})} \\
 &\lesssim& \frac{1}{\epsilon^{\alpha - \gamma}} \sum\limits_{j = 0}^\gamma \left(\frac{\left| a\right|}{\epsilon}\right )^j \left\| \left(\frac{y^1 - a}{\epsilon}\right)^{\gamma - j} \frac{\partial^{\alpha + \beta} F_0}{\partial (y^1)^\alpha \partial R^\beta}(\frac{\cdot - a}{\epsilon};R) \right\|_{L^2(\mathbb{R})} \\
 &\lesssim& \frac{1}{\epsilon^{\alpha - \gamma - \frac{1}{2}}} \left( 1 + \frac{\left| a(y^0)\right|}{\epsilon} \right)^\gamma
\end{eqnarray*}
where we did a change of variables and used the exponential decay of $F_0$ and it's derivatives to obtain the last inequality.
\bigskip

We estimate (\ref{partial F_1 estimate}) in the same way.
\begin{flushright}
 $\Box$
\end{flushright}

We will use lemma \ref{partial F_i estimate lemma} to prove the more useful estimates
\begin{corollary}\label{partial F_i tilde estimate corollary}
 For $\alpha$, $\gamma \in \mathbb{N}\cup \left\{0\right\}$ and $\beta = 0,1,2$, then for $0 \leq y^0 \leq y_*^0$ and for $i=1,2$ we have that
 \begin{equation}\label{partial F_i tilde estimate}
  \left\| (y^1)^\gamma (\frac{\partial}{\partial y^1})^\alpha (\frac{\partial}{\partial y^0})^\beta \tilde{F}_i  \right\|_{L^2(\mathbb{R})} \lesssim \frac{1}{\epsilon^{\alpha - \gamma - 1/2}} (1 + \delta^{\beta 2} \frac{\left| a''\right|}{\epsilon} )(1 + \frac{\left| a\right|}{\epsilon} + \frac{\left| a'\right|}{\epsilon})^{\gamma + \beta}
 \end{equation}
 where $\delta^{ij}$ is the Kronecker-delta and the constant in the estimate depends on $y_*^0$, but not $\epsilon$.
\end{corollary}

\noindent \underline{Proof of corollary \ref{partial F_i tilde estimate corollary}}: We will show (\ref{partial F_i tilde estimate}) for $i=0$. The same arguments can be used to show (\ref{partial F_i tilde estimate}) for $i=1$.
\bigskip

\noindent \underline{$\beta = 0$}: This directly follows from lemma \ref{partial F_i estimate lemma}.
\bigskip

\noindent \underline{$\beta = 1$}: We have that
\begin{gather}
 \frac{\partial}{\partial y^0} \tilde{F}_0 = -\frac{a'}{\epsilon} \frac{\partial F_0}{\partial y^1} + R' \; \frac{\partial F_0}{\partial R} \label{partial y^0 F_0 tilde}
\end{gather}
where $F_0$ and all of its partial derivatives are evaluated at $(\frac{y^1 - a}{\epsilon};R)$. We suppressed the arguments of these quantities for notational convenience. Estimating $(y^1)^\gamma \partial_{y^0} \partial_{y^1}^\alpha \tilde{F}_0$ first we have
\begin{eqnarray*}
 & & \left\| (y^1)^\gamma \frac{\partial}{\partial y^0} (\frac{\partial}{\partial y^1})^\alpha \tilde{F}_0\right\|_{L^2(\mathbb{R})} \\
 &\leq& \left\| (y^1)^\gamma a' (\frac{\partial}{\partial y^1})^{\alpha + 1} \left[ F_0(\frac{\cdot - a}{\epsilon};R)\right] \right\|_{L^2(\mathbb{R})} + \left\| (y^1)^\gamma R' (\frac{\partial}{\partial y^1})^\alpha \frac{\partial}{\partial R} \left[F_0(\frac{\cdot - a}{\epsilon};R)\right]  \right\|_{L^2(\mathbb{R})} \\
 &\lesssim& \frac{1}{\epsilon^{\alpha - \gamma - 1/2}}\left(1 + \frac{\left| a\right|}{\epsilon} + \frac{\left| a'\right|}{\epsilon}\right)^{\gamma + 1}
\end{eqnarray*}
where we used (\ref{partial F_0 estimate}) to obtain the last inequality.
\bigskip

\noindent \underline{$\beta = 2$}: We have that
\begin{align}
 \frac{\partial^2}{\partial (y^0)^2} \tilde{F}_0 &= -\frac{a''}{\epsilon} \frac{\partial F_0}{\partial y^1} + (\frac{a'}{\epsilon})^2 \frac{\partial^2 F_0}{\partial (y^1)^2} + 2 R' \frac{a'}{\epsilon} \frac{\partial^2 F_0}{\partial y^1 \partial R} + (R')^2 \frac{\partial^2 F_0}{\partial R^2}  + R'' \frac{\partial F_0}{\partial R} \label{partial y^0 y^0 F_0 tilde}
\end{align}
where again $F_0$ and all of its partial derivatives are evaluated at $(\frac{y^1 - a}{\epsilon};R)$. Estimating $(y^1)^\gamma \partial_{y^0}^2 \partial_{y^1}^\alpha \tilde{F}_0$ in the same way as we did when finding (\ref{partial F_i tilde estimate}) for $i=0$ and $\beta = 1$, we have that
\begin{eqnarray*}
 \left\| (y^1)^\gamma \partial_{y^0}^2 \partial_{y^1}^\alpha \tilde{F}_0\right\|_{L^2(\mathbb{R})} 
 &\lesssim& \frac{1}{\epsilon^{\alpha - \gamma - 1/2}} (1 + \left| \frac{a''}{\epsilon} \right|) (1 + \frac{\left| a\right|}{\epsilon} + \left| \frac{a'}{\epsilon} \right|)^{\gamma + 2}
\end{eqnarray*}
where we used lemma \ref{partial F_0 estimate} to obtain the last inequality.
\begin{flushright}
 $\Box$
\end{flushright}

To begin the energy estimate, we will use the following divergence identity.
\begin{lemma}\label{divergence identity lemma}
 \begin{equation}\label{divergence identity}
  \partial_{y^0}\xi\cdot \left[ \frac{m^2}{n^2} \partial_{y^0 y^0}\xi + B^\alpha \partial_\alpha \xi + L_\epsilon(\tilde{F}_0,R)\xi \right] = div_{y^0,y^1} \overrightarrow{X} + Y
 \end{equation}
 where $$\overrightarrow{X} = (e(\xi), -\partial_{y^0}\xi \cdot \partial_{y^1}\xi)$$ $$Y = -\frac{1}{\epsilon^2} \xi \cdot \left[ \partial_{y^0} \operatorname{Hess}_\Phi W(\tilde{F}_0,R)\right] \xi + B^\alpha \partial_{y^0}\xi \cdot \partial_\alpha \xi - \frac{1}{2} \partial_{y^0} (\frac{m^2}{n^2}) \partial_{y^0}\xi^2 $$
\end{lemma}
\noindent We omit the proof of lemma \ref{divergence identity lemma} as the proof is a straightforward computation. Using the divergence identity (\ref{divergence identity}) and (\ref{xi equation}), we have
\begin{eqnarray*}
 \partial_{y^0} E 
 &=& \partial_{y^0}\xi \cdot \partial_{y^1}\xi\bigg|_{-y_*^1}^{y_*^1} - \int\limits_{\left| y^1 \right| \leq y_*^1} \partial_{y^0} \xi \cdot\left[S_{-1} + S_0 + N\right] - \int\limits_{\left| y^1 \right| \leq y_*^1} Y
\end{eqnarray*}
Integrating $\partial_{y^0} E$ with respect to $y^0$ once to get
\begin{equation}\label{energy identity}
 E(y^0) - E(0) = \int\limits_0^{y^0} \partial_{y^0}\xi\cdot \partial_{y^1}\xi \bigg|_{-y_*^1}^{y_*^1} -\int\limits_0^{y^0} \int\limits_{\left| y^1 \right| \leq y_*^1} \partial_{y^0}\xi\cdot \left[S_{-1} + S_0 + N\right] - \int\limits_0^{y^0}\int\limits_{\left| y^1 \right| \leq y_*^1} Y
\end{equation}
This energy identity is the main equation we want to estimate.
\bigskip

We will break the analysis up to simplify things. We will consider each term on the right hand side of (\ref{energy identity}) individually, estimate them, and then in the end add all of the individual estimates back up to obtain the desired estimate.

\begin{lemma}\label{S_-1 estimate lemma}
 \begin{equation}\label{S_-1 estimate}
  -\int\limits_0^{y^0} \int\limits_{\left| y^1 \right| \leq y_*^1} \partial_{y^0}\xi \cdot S_{-1} \lesssim (\epsilon E^{1/2} + e^{-\alpha \frac{y_*^1 - a}{\epsilon}}) \sqrt{\epsilon} \underline{A}\bigg|_0^{y^0} + \int\limits_0^{y^0} (\epsilon E^{1/2} + e^{-\alpha \frac{y_*^1 - a}{\epsilon}}) \sqrt{\epsilon} \underline{A}
 \end{equation}
\end{lemma}

\noindent \underline{Proof}: Recall that $\tilde{F}_1$ solves (\ref{translated F_1 equation}) $$\epsilon L_\epsilon(\tilde{F}_0,R)\tilde{F}_1 = H(R) \partial_{y^1} \tilde{F}_0 - (\frac{y^1 - a(y^0)}{\epsilon}) m(y^0) \partial_R w(\tilde{F}_0,R) $$ Using (\ref{B^1}) and (\ref{mean curvature}) we see that $$B^1 = -H(R) + O(y^1)$$ Recall the definition of $S_{-1}$ (\ref{S_-1})
\begin{equation}\label{S_-1 mean curvature}
 S_{-1} = \left[ H(R) + B^1\right]\partial_{y^1}\tilde{F}_0 + \frac{a}{\epsilon} m(y^0) \partial_R w(\tilde{F}_0,R)
\end{equation}
Integrating by parts in $y^0$ we have 
\begin{equation}\label{S_-1 integrate by parts}
 -\int\limits_0^{y^0}\int\limits_{\left| y^1 \right| \leq y_*^1} \partial_{y^0}\xi \cdot S_{-1} = - \int\limits_{\left| y^1 \right| \leq y_*^1}  \xi \cdot  S_{-1}\bigg|_0^{y^0} + \int\limits_0^{y^0}\int\limits_{\left| y^1 \right| \leq y_*^1} \xi \cdot \partial_{y^0} S_{-1}
\end{equation}
\bigskip

For $j=0,1$ we will need to estimate
\begin{equation}\label{S_-1 cauchy schwarz}
 \int\limits_{-y_*^1}^{y_*^1} \xi \cdot \partial_{y^0}^j S_{-1}
 \lesssim \left\|\xi\right\|_{L^2(-y_*^1,y_*^1)}  \left\| \partial_{y^0}^j S_{-1}\right\|_{L^2(\mathbb{R})} 
 \lesssim (\epsilon E^{1/2} + e^{-\alpha \frac{y_*^1 - a}{\epsilon}})\left\| \partial_{y^0}^j S_{-1}\right\|_{L^2(\mathbb{R})}
\end{equation}
We estimate $\left\| S_{-1}\right\|_{L^2(\mathbb{R})}$ first as
\begin{gather}\label{nonderivative S_-1 estimate}
 \left\| S_{-1}\right\|_{L^2(\mathbb{R})}
 \lesssim \left\| (y^1) \partial_{y^1} \tilde{F}_0 \right\|_{L^2(\mathbb{R})} + \left\| \frac{a(y^0)}{\epsilon} m(y^0) \partial_R w(\tilde{F}_0,R)\right\|_{L^2(\mathbb{R})}
 \lesssim \sqrt{\epsilon} \underline{A}
\end{gather}
where we used corollary \ref{partial F_i tilde estimate corollary} to obtain the last inequality. We estimate the second term of (\ref{S_-1 integrate by parts}) as
\begin{gather}\label{derivative S_-1 estimate}
 \left\| \partial_{y^0} S_{-1}\right\|_{L^2(\mathbb{R})}
 \lesssim \left\| \partial_{y^0}\left[ \left( H(R) + B^1 \right)\partial_{y^1}\tilde{F}_0 \right]\right\|_{L^2(-y_*^1,y_*^1)} + \left\| \partial_{y^0} \left[ \frac{a}{\epsilon} m(y^0) \partial_R w(\tilde{F}_0,R) \right] \right\|_{L^2(-y_*^1,y_*^1)} 
 \lesssim \sqrt{\epsilon} \underline{A}
\end{gather}
where we again used corollary \ref{partial F_i tilde estimate corollary} to obtain the last inequality. Combining (\ref{S_-1 cauchy schwarz}), (\ref{nonderivative S_-1 estimate}), and (\ref{derivative S_-1 estimate}) finishes the proof.
\begin{flushright}
 $\Box$
\end{flushright}

\begin{lemma}\label{S_0 estimate lemma}
 \begin{equation}\label{S_0 estimate}
  -\int\limits_0^{y^0}\int\limits_{\left| y^1 \right| \leq y_*^1} \partial_{y^0}\xi\cdot S_0 \lesssim (\epsilon E^{1/2} + e^{-\alpha \frac{y^1 - a}{\epsilon}}) \sqrt{\epsilon} \underline{A}\bigg|_0^{y^0} + \int\limits_0^{y^0} (1 + \frac{\left| a'' \right|}{\epsilon}) (\epsilon E^{1/2} + e^{-\alpha \frac{y^1 - a}{\epsilon}}) \sqrt{\epsilon} \underline{A}
 \end{equation}
\end{lemma}
\bigskip

\noindent \underline{Proof}: Using the definition of $S_0$, see (\ref{S_0}), the left hand side of (\ref{S_0 estimate}) is $$-\int\limits_0^{y^0} \int\limits_{\left| y^1 \right| \leq y_*^1} \partial_{y^0}\xi\cdot \left[ \frac{m^2}{n^2} \partial_{y^0 y^0}(\tilde{F}_0 +\epsilon \tilde{F}_1) + B^0\partial_{y^0}(\tilde{F}_0 + \epsilon \tilde{F}_1) + \epsilon B^1 \partial_{y^1} \tilde{F}_1 \right]$$ We would like to integrate by parts in $y^0$ to move the derivative from $\xi$ to $S_0$ and use corollary \ref{partial F_i tilde estimate corollary}. However, then we'd have to estimate $\partial_{y^0 y^0 y^0}(\tilde{F}_0 + \epsilon \tilde{F}_1)$ which will give rise to $a'''$ terms which we'd rather avoid. So, we take special care when estimating these two problematic terms and proceed as we would like to for the other terms.
\begin{enumerate}
 \item \underline{$\partial_{y^0 y^0} \tilde{F}_0$ term}: Recall that $\tilde{F}_0 = F_0(\frac{y^1 - a(y^0)}{\epsilon};R(y^0))$ and so 
 \begin{eqnarray}\label{partial y^0 y^0 F_0}
  \partial_{y^0 y^0} \tilde{F}_0
  &=& \frac{a''}{\epsilon} \partial_{y^1} F_0 - (\frac{a'}{\epsilon})^2 \partial_{y^1 y^1} F_0 + 2 R' \frac{a'}{\epsilon} \partial_R \partial_{y^1} F_0 - R'' \partial_R F_0 - (R')^2 \partial_{RR} F_0
 \end{eqnarray}
where we use the notation $\partial_{y^1}^\alpha \partial_R^\beta F_0 = \partial_{y^1}^\alpha \partial_R^\beta F_0(\frac{y^1 - a}{\epsilon};R)$ to simplify things. To control the $a''$ term, we use the fact that $\xi\perp \partial_{y^1} F_0$. Differentiating $$ 0 = \int\limits_\mathbb{R} \xi \cdot \partial_{y^1} \tilde{F}_0 $$ with respect to $y^0$ once yields
\begin{equation}\label{partial xi dot partial F_0}
  \int\limits_{\mathbb{R}} \partial_{y^0} \xi \cdot \partial_{y^1} F_0 = \frac{a'}{\epsilon} \int\limits_{\mathbb{R}} \xi \cdot \partial_{y^1 y^1} F_0 - R' \int\limits_{\mathbb{R}} \xi \cdot \partial_R\partial_{y^1} F_0
\end{equation}
To use this we first recall that $m = (1 - (R')^2)^{-1/2}$ and $n = 1 + y^1 m^3 R''$. Thus,
$$\int\limits_0^{y^0} \int\limits_{\left| y^1 \right| \leq y_*^1}  \frac{m^2}{n^2} \frac{a''}{\epsilon} \partial_{y^0} \xi \cdot \partial_{y^1} F_0 = \int\limits_0^{y^0} \int\limits_{\left| y^1 \right| \leq y_*^1} m^2 \frac{a''}{\epsilon} \partial_{y^0} \xi \cdot \partial_{y^1} F_0 + \int\limits_0^{y^0} \int\limits_{\left| y^1 \right| \leq y_*^1} \left[ \frac{m^2}{n^2} - m^2\right] \frac{a''}{\epsilon} \partial_{y^0} \xi \cdot \partial_{y^1} F_0$$ Using (\ref{partial xi dot partial F_0}) we control the first term, the $m^2$ term, as follows
\begin{eqnarray*}
 \left| \int\limits_0^{y^0} \int\limits_{\left| y^1 \right| \leq y_*^1} m^2 \frac{a''}{\epsilon} \partial_{y^0}\xi \cdot \partial_{y^1} F_0 \right|
 &\lesssim& \int\limits_0^{y^0} \left| m\right|^2 \frac{\left| a''\right|}{\epsilon} \left| \; \int\limits_{\left| y^1 \right| \leq y_*^1} \partial_{y^0}\xi\cdot \partial_{y^1} F_0\right| \\
 &\lesssim& \int\limits_0^{y^0} \left| \frac{a' a''}{\epsilon^2} \right| \int\limits_{\mathbb{R}} \left| \xi \cdot \partial_{y^1 y^1} F_0\right| + \int\limits_0^{y^0} \left| \frac{a''}{\epsilon} \right|\int\limits_{\mathbb{R}} \left| \xi\cdot \partial_R \partial_{y^1} F_0\right|
\end{eqnarray*}
where we used the fact that $\left| m\right|^2 \lesssim 1$ for all $y^0 \leq y_*^0$ to go from the 1st line to the 2nd. Remember that we are using the notation $\partial_{y^1}^\alpha \partial_{R}^\beta F_0 = \partial_{y^1}^\alpha \partial_R^\beta F_0(\frac{y^1 - a}{\epsilon};R)$. We then use the Cauchy-Schwarz inequality, corollary \ref{partial F_i tilde estimate corollary}, lemma \ref{asymptotics}, and (\ref{energy-spectral estimate}) to conclude that
\begin{eqnarray*}
 \left| \int\limits_0^{y^0} \int\limits_{\left| y^1 \right| \leq y_*^1} m^2 \frac{a''}{\epsilon} \partial_{y^0} \xi\cdot \partial_{y^1} F_0\right| 
 &\lesssim& \int\limits_0^{y^0} \frac{\left| a'' \right|}{\epsilon} (\epsilon E^{1/2} + e^{-\alpha \frac{y^1 - a}{\epsilon}}) \sqrt{\epsilon} \underline{A}
\end{eqnarray*}
To control the second term (i.e. the $\frac{m^2}{n^2} - m^2$ term), observe that $\frac{m^2}{n^2} - m^2 = O(y^1)$. We use the Cauchy-Schwarz inequality and the definition of energy to get
\begin{eqnarray*}
 \left| \int\limits_0^{y^0} \int\limits_{\left| y^1 \right| \leq y_*^1} \left[ \frac{m^2}{n^2} - m^2 \right] \frac{a''}{\epsilon} \partial_{y^0}\xi\cdot \partial_{y^1} F_0\right|
 &\lesssim& \int\limits_0^{y^0} \frac{\left| a'' \right|}{\epsilon} \left\| \partial_{y^0} \xi\right\|_{L^2(-y_*^1,y_*^1)} \left\| y^1 \partial_{y^1} F_0 \right\|_{L^2(\mathbb{R})} \\
 &\lesssim& \int\limits_0^{y^0} \frac{\left| a'' \right|}{\epsilon} \epsilon^{3/2} E^{1/2} \underline{A} 
\end{eqnarray*}

To summarize, we have the following estimate
\begin{eqnarray*}
 \int\limits_0^{y^0} \int\limits_{\left| y^1 \right| \leq y_*^1} \frac{m^2}{n^2} \frac{a''}{\epsilon} \partial_{y^0} \xi\cdot\partial_{y^1} F_0 
 &\lesssim& \int\limits_0^{y^0} \frac{\left| a'' \right|}{\epsilon} (\epsilon E^{1/2} + e^{-\alpha \frac{y^1 - a}{\epsilon}}) \sqrt{\epsilon} \underline{A}
\end{eqnarray*}

We deal with the rest of the terms of (\ref{partial y^0 y^0 F_0}) by shifting the $\partial_{y^0}$ off of $\xi$ and use (\ref{energy-spectral estimate}) along with (\ref{partial F_i estimate lemma}). We estimate the $\partial_{y^1 y^1} F_0$ term of (\ref{partial y^0 y^0 F_0}) and the estimation of the other three terms of (\ref{partial y^0 y^0 F_0}) are done in the same way yielding the same bounds. To this end, we estimate the $\partial_{y^1 y^1} F_0$ term as
\begin{eqnarray*}
 & &\int\limits_0^{y^0} \int\limits_{\left| y^1 \right| \leq y_*^1} \frac{m^2}{n^2} \partial_{y^0}\xi\cdot\left[ -(\frac{a'}{\epsilon})^2 \partial_{y^1 y^1}F_0 \right] \\
 &=& - \int\limits_{\left| y^1 \right| \leq y_*^1} \frac{m^2}{n^2}(\frac{a'}{\epsilon})^2 \xi\cdot \partial_{y^1 y^1} F_0\bigg|_0^{y^0} + \int\limits_0^{y^0}\int\limits_{\left| y^1 \right| \leq y_*^1} \xi\cdot \partial_{y^0}\left[ \frac{m^2}{n^2} (\frac{a'}{\epsilon})^2 \partial_{y^1 y^1} F_0 \right] \\
 &\lesssim& (\frac{a'}{\epsilon})^2 \left\|\xi \right\|_{L^2(-y_*^1,y_*^1)} \left\| \partial_{y^1 y^1} F_0\right\|_{L^2(\mathbb{R})}\bigg|_0^{y^0} + \int\limits_0^{y^0} \left\|\xi\right\|_{L^2(-y_*^1,y_*^1)}\left\| \partial_{y^0}\left[ \frac{m^2}{n^2} (\frac{a'}{\epsilon})^2 \partial_{y^1 y^1} F_0 \right] \right\|_{L^2(\mathbb{R})} \\
 &\lesssim& (\epsilon E^{1/2} + e^{-\alpha \frac{y^1 - a}{\epsilon}}) \sqrt{\epsilon} \underline{A} \bigg|_0^{y^0} + \int\limits_0^{y^0} (1 + \frac{\left|a''\right|}{\epsilon})(\epsilon E^{1/2} + e^{-\alpha \frac{y^1 - a}{\epsilon}}) \sqrt{\epsilon} \underline{A}
\end{eqnarray*}

Putting together the above estimates yields
\begin{equation}\label{partial y^0 y^0 F_0 estimate}
 - \int\limits_0^{y^0} \int\limits_{\left| y^1 \right| \leq y_*^1} \frac{m^2}{n^2} \partial_{y^0}\xi\cdot \partial_{y^0 y^0} \tilde{F}_0 \lesssim (\epsilon E^{1/2} + e^{-\alpha \frac{y^1 - a}{\epsilon}}) \sqrt{\epsilon} \underline{A} \bigg|_0^{y^0} + \int\limits_0^{y^0} (1 + \frac{\left|a''\right|}{\epsilon})(\epsilon E^{1/2} + e^{-\alpha \frac{y^1 - a}{\epsilon}}) \sqrt{\epsilon} \underline{A}
\end{equation}

\item \underline{$\partial_{y^0 y^0} \tilde{F}_1$ term}: Since $\frac{m^2}{n^2} \lesssim 1$ on $(0,y_*^0)\times (-y_*^1,y_*^1)$, then after using Cauchy-Schwarz and corollary \ref{partial F_i tilde estimate corollary} we have
\begin{gather}\label{partial y_0 y_0 F_1 estimate}
 -\epsilon \int\limits_0^{y^0} \int\limits_{\left| y^1 \right| \leq y_*^1} \frac{m^2}{n^2} \partial_{y^0}\xi\cdot \partial_{y^0 y^0}\tilde{F}_1
 = \epsilon \int\limits_0^{y^0} E^{1/2} \left\| \partial_{y^0 y^0} \tilde{F}_1\right\|_{L^2(\mathbb{R})}
 \lesssim \int\limits_0^{y^0} (1 + \frac{\left| a'' \right|}{\epsilon}) \epsilon^{3/2} E^{1/2} \underline{A} 
\end{gather}

\item \underline{$B^0 \partial_{y^0} \tilde{F}_0$}: Using the boundedness of $B^0$ on $(0,y_*^0)\times (-y_*^1,y_*^1)$ (defined in (\ref{B^0})), we have
\begin{eqnarray}\label{B^0 partial y_0 F_0 estimate}
 & &-\int\limits_0^{y^0} \int\limits_{\left| y^1 \right|\leq y_*^1} B^0\partial_{y^0}\xi\cdot \partial_{y^0}\tilde{F}_0 \nonumber \\
 &=&  -\int\limits_{\left| y^1 \right|\leq y_*^1} B^0 \xi\cdot \partial_{y^0} \tilde{F}_0\bigg|_0^{y^0} + \int\limits_0^{y^0}\int\limits_{\left| y^1 \right|\leq y_*^1} \xi\cdot \partial_{y^0} \left[B^0 \partial_{y^0}\tilde{F}_0\right] \nonumber \\
 &\lesssim& (\epsilon E^{1/2} + e^{-\alpha \frac{y^1 - a}{\epsilon}})\sqrt{\epsilon}\underline{A}\bigg|_0^{y^0} + \int\limits_0^{y^0} (1 + \frac{\left|a''\right|}{\epsilon})  (\epsilon E^{1/2} + e^{-\alpha \frac{y^1-a}{\epsilon}}) \sqrt{\epsilon}\underline{A} 
\end{eqnarray}
where we used corollary \ref{partial F_i tilde estimate corollary} to obtain the last inequality.

\item \underline{$B^0\partial_{y^0}\tilde{F}_1$ term}: We estimate this term in the same way as the $B^0 \partial_{y^0} \tilde{F}_0$ term. That is,
\begin{eqnarray}\label{B^0 partial y_0 F_1 estimate}
 -\epsilon\int\limits_0^{y^0}\int\limits_{-y_*^1}^{y_*^1} B^0 \partial_{y^0}\xi\cdot\partial_{y^0}\tilde{F}_1 &\lesssim& \int\limits_0^{y^0} (\epsilon E^{1/2} + e^{-\alpha \frac{y^1 - a}{\epsilon}})\sqrt{\epsilon}\underline{A}
\end{eqnarray}

\item \underline{$B^1 \partial_{y^1} \tilde{F}_1$ term}: Using the boundedness of $B^1$ on $(0,y_*^0)\times (-y_*^1,y_*^1)$ we get
\begin{eqnarray}\label{B^1 partial y_1 F_1 estimate}
 & &-\int\limits_0^{y^0}\int\limits_{\left| y^1 \right|\leq y_*^1} \epsilon B^1 \partial_{y^0}\xi\cdot\partial_{y^1}\tilde{F}_1 \nonumber \\
 &=&-\int\limits_{\left| y^1 \right|\leq y_*^1} \xi\cdot\left[\epsilon B^1\partial_{y^1} \tilde{F}_1\right]\bigg|_0^{y^0} + \int\limits_0^{y^0}\int\limits_{\left| y^1 \right|\leq y_*^1}\xi\cdot\partial_{y^0}\left[\epsilon B^1\partial_{y^1} \tilde{F}_1\right] \nonumber \\
 &\lesssim& (\epsilon E^{1/2} + e^{-\alpha \frac{y^1-a}{\epsilon}})\sqrt{\epsilon}\underline{A} \bigg|_0^{y^0} + \int\limits_0^{y^0}(\epsilon E^{1/2} + e^{-\alpha \frac{y^1-a}{\epsilon}})\sqrt{\epsilon}\underline{A}
\end{eqnarray}
where we used corollary \ref{partial F_i tilde estimate corollary} to obtain the last inequality.
\end{enumerate}
\bigskip

Putting together the estimates obtained from steps 1-5 we get the desired estimate.
\begin{flushright}
 $\Box$
\end{flushright}

\begin{lemma}\label{N estimate lemma}
 \begin{equation}\label{N estimate}
  \left| \int\limits_0^{y^0}\int\limits_{\left| y^1 \right| \leq y_*^1} \partial_{y^0}\xi\cdot N\right|
 \lesssim \left[\frac{y_*^1}{(R - y_*^1 m)^3} E + \left(\sqrt{\epsilon}^{3}\underline{A} + \frac{1}{\epsilon} e^{-\alpha \frac{y_*^1 - a}{\epsilon}}\right) \sqrt{E} + \sqrt{\epsilon} \sqrt{E}^{3}\right]_0^{y^0} + \int\limits_0^{y^0} \epsilon \sqrt{E}
 \end{equation}
\end{lemma}

Note that we have a factor $\frac{y_*^1}{(R - y_*^1 m)^3}$ in the above estimate. We'll need to pick $y_*^1$ sufficiently small so that the constant coming from this estimate multiplied by these two factors is less than 1. We need this to be able to close the bootstrap argument.
\bigskip

\noindent \underline{Proof}: Recalling (\ref{N}), we have $$N = -\frac{1}{\epsilon^2}\left[w(\tilde{F}_0 + \tilde{F}_\xi, R+y^1 m) - w(\tilde{F}_0,R) - \operatorname{Hess}_\Phi(\tilde{F}_0,R)\tilde{F}_\xi - y^1 m(y^0) \partial_R w(\tilde{F}_0,R) )\right]$$ where, recall, $\tilde{F}_\xi = \epsilon \tilde{F}_1 + \xi$. Using the identity $g(t) = g(0) + g'(0)t +\int\limits_0^1 (1-t)g''(t)dt$ we can rewrite $N$ as $$N=-\frac{1}{\epsilon^2}\int\limits_0^1 (1-t)\frac{d^2}{dt^2} w(\tilde{F}_0 + t\tilde{F}_\xi, R + t y^1 m) dt$$ Examining (\ref{N estimate}) and integrating by parts with respect to $y^0$ one has that $$-\int\limits_0^{y^0}\int\limits_{-y_*^1}^{y_*^1}\partial_{y^0}\xi\cdot N = -\int\limits_{-y_*^1}^{y_*^1} \xi\cdot N\bigg|_0^{y^0} + \int\limits_0^{y^0}\int\limits_{-y_*^1}^{y_*^1}\xi\cdot\partial_{y^0} N$$ Using Cauchy-Schwarz, we have 
\begin{equation}\label{integrate by parts cauchy-schwarz N}
 -\int\limits_0^{y^0}\int\limits_{-y_*^1}^{y_*^1}\partial_{y^0}\xi\cdot N \lesssim \epsilon E^{1/2}\left\| N\right\|_{L^2(-y_*^1,y_*^1)} \bigg|_0^{y^0} + \int\limits_0^{y^0} \epsilon E^{1/2} \left\|\partial_{y^0} N\right\|_{L^2(-y_*^1,y_*^1)}
\end{equation}
We are then left to estimate $\left\| N\right\|_{L^2(-y_*^1,y_*^1)}$ and $\left\|\partial_{y^0} N\right\|_{L^2(-y_*^1,y_*^1)}$.
\bigskip

Before that, we examine the $\frac{d^2}{dt^2}$ term in $N$. Recall that $\tilde{F}_\xi(y^0,y^1)$ is a two-component vector. For $\tilde{F}_\xi = \left( (\tilde{F}_\xi)_\phi, (\tilde{F}_\xi)_\sigma \right)$, we have that
\begin{eqnarray*}
 & &\frac{d^2}{dt^2} w(\tilde{F}_0 + t\tilde{F}_\xi,R + t y^1 m) \nonumber \\
 &=& \left( (\tilde{F}_\xi)_\phi \partial_\phi + (\tilde{F}_\xi)_\sigma \partial_\sigma\right) \operatorname{Hess}_\Phi W(\tilde{F}_0 + t\tilde{F}_\xi, R + t y^1 m)\tilde{F}_\xi \nonumber \\
 & & + 2 y^1 m(y^0) \partial_R \operatorname{Hess}_\Phi W(\tilde{F}_0 + t\tilde{F}_\xi, R + t y^1 m)\tilde{F}_\xi + (y^1 m)^2 \partial_{RR} w(\tilde{F}_0 + t\tilde{F}_\xi, R + t y^1 m) \nonumber
\end{eqnarray*}

\noindent \underline{$\left\| \partial_{y^0} N\right\|_{L^2(-y_*^1,y_*^1)}$ Estimate}: For the purpose of our result we only need $\left\| \partial_{y^0} N\right\|_{L^2(-y_*^1,y_*^1)} \lesssim 1$. This is straightforward as $\tilde{F}_0$ and $\tilde{F}_\xi$, and all time derivatives of these quantities have bounded $L^2$-norm.
\bigskip

\noindent \underline{$\left\| N\right\|_{L^2(-y_*^1,y_*^1)}$ Estimate}:
\begin{enumerate}
 \item \underline{$(\tilde{F}_\xi)_\phi$ term}:
 \begin{eqnarray*}
  \left\|\int\limits_0^1 (1-t)(\tilde{F}_\xi)_\phi \partial_\phi \operatorname{Hess}_\Phi W(\tilde{F}_0 + t \tilde{F}_\xi,R+t y^1 m)\tilde{F}_\xi dt\right\|_{L^2(-y_*^1,y_*^1)}
  &\lesssim& \left\| (\tilde{F}_\xi)_\phi\right\|_\infty\left\| \tilde{F}_\xi\right\|_{L^2(-y_*^1,y_*^1)}
 \end{eqnarray*}
 where we used the boundedness of $\tilde{F}_\xi$ to control the operator norm $\left\| \partial_\phi \operatorname{Hess}_\Phi W\right\|$. Next, we use Gagliardo-Nirenberg to show that
 \begin{eqnarray*}
  \left\| \tilde{F}_\xi\right\|_\infty 
  &\lesssim& \epsilon \left\| \tilde{F}_1\right\|_\infty + \left\|\xi\right\|_\infty \\
  &\lesssim& \epsilon + \left\|\xi\right\|_{L^2(-y_*^1,y_*^1)}^{1/2} \left\|\partial_{y^1} \xi\right\|_{L^2(-y_*^1,y_*^1)}^{1/2} \\
  &\lesssim& \epsilon + \epsilon^{1/2} E^{1/2}
 \end{eqnarray*}
and we estimate
\begin{eqnarray*}
 \left\|\tilde{F}_\xi\right\|_{L^2(-y_*^1,y_*^1)} 
 &\lesssim& \epsilon \left\|\tilde{F}_1\right\|_{L^2(\mathbb{R})} + \left\|\xi\right\|_{L^2(-y_*^1,y_*^1)} \\
 &\lesssim& \epsilon^{3/2} + (\epsilon E^{1/2} + e^{-\alpha \frac{y_*^1 - a}{\epsilon}})
\end{eqnarray*}
Putting these together, we obtain the estimate
\begin{eqnarray*}
 & & \left\| -\frac{1}{\epsilon^2}\int\limits_0^1 (1-t)(\tilde{F}_\xi)_\phi \partial_\phi \operatorname{Hess}_\Phi W(\tilde{F}_0 + t \tilde{F}_\xi,R+t y^1 m)\tilde{F}_\xi dt\right\|_{L^2(-y_*^1,y_*^1)} \\
 &\lesssim& \epsilon^{1/2} + (1 + \frac{1}{\epsilon^{3/2}}e^{-\alpha \frac{y_*^1 - a}{\epsilon}}) E^{1/2} + \frac{1}{\sqrt{\epsilon}} E + \frac{1}{\epsilon^2} e^{-\alpha \frac{y_*^1 - a}{\epsilon}}
\end{eqnarray*}

\item \underline{$(\tilde{F}_\xi)_\sigma$ term}: We estimate this term in the same way we did the first. Thus,
\begin{eqnarray*}
 & & \left\| -\frac{1}{\epsilon^2}\int\limits_0^1 (1-t)(\tilde{F}_\xi)_\sigma \partial_\sigma \operatorname{Hess}_\Phi W(\tilde{F}_0 + t \tilde{F}_\xi,R+t y^1 m)\tilde{F}_\xi dt\right\|_{L^2(-y_*^1,y_*^1)} \\
 &\lesssim& \epsilon^{1/2} + (1 + \frac{1}{\epsilon^{3/2}}e^{-\alpha \frac{y_*^1 - a}{\epsilon}}) E^{1/2} + \frac{1}{\sqrt{\epsilon}} E + \frac{1}{\epsilon^2} e^{-\alpha \frac{y_*^1 - a}{\epsilon}} 
\end{eqnarray*}

\item \underline{$2 y^1 m(y^0) \partial_R \operatorname{Hess}_\Phi W \tilde{F}_\xi$ term}: 
\begin{eqnarray*}
 & & \left\| -\frac{1}{\epsilon^2} \int\limits_0^1 (1-t)2 y^1 m \partial_R \operatorname{Hess}_\Phi W(\tilde{F}_0 + t \tilde{F}_\xi,R+t y^1 m) \tilde{F}_\xi\right\|_{L^2(-y_*^1,y_*^1)} \\
 &\lesssim& \frac{1}{\epsilon^2} \left\| \frac{y^1}{(R - y_*^1 m)^3} \tilde{F}_\xi\right\|_{L^2(-y_*^1,y_*^1)} \\
 &\lesssim&\frac{1}{\epsilon^2}\left( \epsilon^{5/2}\underline{A} + \epsilon \frac{y_*^1}{(R - y_*^1 m)^3} E^{1/2} + e^{-\alpha \frac{y_*^1-a}{\epsilon}} \right)
\end{eqnarray*}
where we made use of corollary \ref{partial F_i tilde estimate corollary} to obtain the estimate.

\item \underline{$(y^1 m(y^0))^2 \partial_{RR}w$ term}: Since $$\partial_{RR} w(\tilde{F}_0 + t \tilde{F}_\xi,R + t y^1 m) = \left( \begin{array}{c} 0 \\ \frac{3 d^2}{(R + t y^1 m)^4}(\tilde{s}_0 + t(\tilde{F}_\xi)_\sigma) \end{array}\right)$$ then we have
\begin{eqnarray*}
 & &\left\| -\frac{1}{\epsilon^2} \int\limits_0^1(1-t)(y^1 m)^2 \partial_{RR} w(\tilde{F}_0 + t\tilde{F}_\xi,R + t y^1 m)\right\|_{L^2(-y_*^1,y_*^1)} \\
 &\lesssim& \frac{1}{\epsilon^2} \left(\left\| (y^1)^2 \tilde{s}_0\right\|_{L^2(\mathbb{R})} + \epsilon \left\| (y^1)^2 \tilde{s}_1\right\|_{L^2(\mathbb{R})} + \left\| \frac{(y^1)^2}{(R - y^1 m)^4} \xi\right\|_{L^2(-y_*^1,y_*^1)} \right) \\
 &\lesssim& \frac{1}{\epsilon^2}\left(\epsilon^{5/2} \underline{A} + \epsilon^{7/2}\underline{A} + \epsilon \frac{(y_*^1)}{(R - y_*^1 m)^3}  E^{1/2} + e^{-\alpha \frac{y_*^1-a}{\epsilon}} \right)
\end{eqnarray*}
\end{enumerate}
where again we made use of corollary \ref{partial F_i tilde estimate corollary} to obtain the estimate.
\bigskip

Putting together estimate (\ref{integrate by parts cauchy-schwarz N}) and the estimates from steps 1-4, we obtain (\ref{N estimate}).
\begin{flushright}
 $\Box$
\end{flushright}

\begin{lemma}\label{Y Estimate Lemma}
 \begin{equation}\label{Y Estimate}
  -\int\limits_0^{y^0} \int\limits_{-y_*^1}^{y_*^1} Y \lesssim \int\limits_0^{y^0}\left( E + \frac{1}{\epsilon^2} e^{-\alpha \frac{y^1 - a}{\epsilon}}\right) 
 \end{equation}
\end{lemma}

\noindent \underline{Proof}: (See lemma \ref{divergence identity} for the definition of $Y$)
\begin{eqnarray*}
 -\int\limits_0^{y^0} \int\limits_{-y_*^1}^{y_*^1} Y
 &=& \frac{1}{\epsilon^2} \int\limits_0^{y^0}\int\limits_{-y_*^1}^{y_*^1} \xi\cdot \left[\partial_{y^0}\operatorname{Hess}_\Phi W(\tilde{F}_0,R)\right] \xi - \int\limits_0^{y^0}\int\limits_{\left|y^1\right|\leq y_*^1} B^\alpha\partial_{y^0}\xi\cdot\partial_\alpha \xi - \int\limits_0^{y^0} \int\limits_{-y_*^1}^{y_*^1} \frac{1}{2}\partial_{y^0}\left(\frac{m^2}{n^2}\right) \partial_{y^0}\xi^2 \\
 &\lesssim& \int\limits_0^{y^0}\left( E + \frac{1}{\epsilon^2} e^{-\alpha \frac{y^1 - a}{\epsilon}}\right)
\end{eqnarray*}
where we used the boundedness of the operator $\partial_{y^0} \operatorname{Hess}_\Phi W(\tilde{F}_0,R)$ and the $B^\alpha$'s on $(0,y_*^0)\times (-y_*^1,y_*^1)$ and (\ref{energy-spectral estimate}) to obtain the estimate.
\begin{flushright}
 $\Box$
\end{flushright}

\begin{lemma}\label{partial y^0 xi partial y^1 xi estimate lemma}
 \begin{equation}\label{partial y^0 xi partial y^1 xi estimate}
  \int\limits_0^{y^0} \partial_{y^0} \xi(y^0,\pm y_*^1)\cdot\partial_{y^1}\xi(y^0,\pm y_*^1) \lesssim \int\limits_0^{y^0} \frac{1}{\epsilon} e^{-\alpha \frac{y^1 - a}{\epsilon}} \underline{A}
 \end{equation}
\end{lemma}

\noindent \underline{Proof}: Recall that $$\xi(y^0,\pm y_*^1) = \left(\begin{array}{c} \pm 1 \\ 0 \end{array}\right) - \tilde{F}_0(y^0,\pm y_*^1) - \epsilon \tilde{F}_1(y^0,\pm y_*^1)$$ Differentiating $\tilde{F}_0$ and $\tilde{F}_1$ with respect to $y^0$ and using (\ref{asymptotics}), then
\[
 \left| \partial_{y^0}\xi(y^0,\pm y_*^1) \right| \lesssim e^{-\alpha \frac{y_*^1 - a(y^0)}{\epsilon}} \underline{A}
\]
Using (\ref{asymptotics}), we also have that
\begin{gather*}
 \left| \partial_{y^1} \xi(y^0,\pm y_*^1)\right|
 = \left| -\frac{1}{\epsilon} \partial_{y^1} F_0(\frac{y_*^1 - a(y^0)}{\epsilon};R(y^0)) - \partial_{y^1} F_1(\frac{y_*^1 - a(y^0)}{\epsilon};R(y^0),R'(y^0))\right|
 \lesssim \frac{1}{\epsilon} e^{-\alpha\frac{y_*^1 - a(y^0)}{\epsilon}}
\end{gather*}
\begin{flushright}
 $\Box$
\end{flushright}

Combining the estimates obtained from lemma \ref{S_-1 estimate lemma} to lemma \ref{partial y^0 xi partial y^1 xi estimate lemma} allows us to conclude the proof of theorem \ref{energy estimate theorem}.

\subsection{Proof of Bounded Shift Theorem (Theorem \ref{bounded shift theorem})}
To prove this we will use the fact that $\xi\perp \partial_{y^1} \tilde{F}_0$. Differentiate this quantity with respect to $y^0$ twice to get
\begin{eqnarray}\label{xi perp partial y^1 F_0}
 0
 &=& \partial_{y^0 y^0} \int\limits_{\mathbb{R}} \xi\cdot\partial_{y^1}\tilde{F}_0 \nonumber \\
 &=& \int\limits_{\mathbb{R}} \partial_{y^0 y^0}\xi\cdot \partial_{y^1}\tilde{F}_0 + 2 \int\limits_{\mathbb{R}} \partial_{y^0}\xi\cdot\partial_{y^0 y^1}\tilde{F}_0 + \int\limits_{\mathbb{R}} \xi\cdot \partial_{y^0 y^0}\partial_{y^1}\tilde{F}_0 \nonumber \\
 &=& \int\limits_{\mathbb{R}}\left( \partial_{y^0 y^0} \xi\cdot\partial_{y^1}\tilde{F}_0 + 2 \partial_{y^0}\xi\cdot \partial_{y^0}\partial_{y^1}\tilde{F}_0 - \partial_{y^1}\xi\cdot\partial_{y^0 y^0} \tilde{F}_0 \right) \nonumber \\
 \int\limits_{\left| y^1\right| \leq y_*^1} \partial_{y^0 y^0} \xi\cdot \partial_{y^1}\tilde{F}_0 &=& \int\limits_\mathbb{R} \left( 2 \partial_{y^0}\xi\cdot \partial_{y^0}\partial_{y^1}\tilde{F}_0 - \partial_{y^1}\xi\cdot\partial_{y^0 y^0} \tilde{F}_0 \right) - \int\limits_{\left| y^1\right| > y_*^1} \partial_{y^0 y^0}\xi \cdot\partial_{y^1} \tilde{F}_0
\end{eqnarray}
where we integrated by parts to move the $\partial_{y^1}$ off of the $\partial_{y^0 y^0} \partial_{y^1} \tilde{F}_0$ term onto the $\xi$ term to obtain the second last equality. On the other hand, we can use the equation for $\xi$ (\ref{xi equation}) to rewrite the left hand side of (\ref{xi perp partial y^1 F_0}) as
\begin{equation}\label{partial y^0 y^0 xi partial y^1 F_0}
 \int\limits_{\left| y^1\right|\leq y_*^1} \partial_{y^0 y^0}\xi\cdot \partial_{y^1} \tilde{F}_0 = -\int\limits_{-y_*^1}^{y_*^1} \frac{n^2}{m^2}\left[B^\alpha\partial_\alpha\xi + L_\epsilon(\tilde{F}_0,R)\xi + S_{-1} + S_0 + N\right]\cdot \partial_{y^1}\tilde{F}_0
\end{equation}
Examining the $S_0$ term on the right hand side of (\ref{partial y^0 y^0 xi partial y^1 F_0}) more closely, we see that
\begin{eqnarray}\label{S_0 partial y^1 F_0 term}
 \int\limits_{-y_*^1}^{y_*^1} \frac{n^2}{m^2} S_0\cdot\partial_{y^1} \tilde{F}_0
 &=& \int\limits_{-y_*^1}^{y_*^1}\left( \partial_{y^0 y^0}\tilde{F}_0 + \epsilon \partial_{y^0 y^0}\tilde{F}_1 + \frac{n^2}{m^2} \partial_{y^0}(\tilde{F}_0 + \epsilon \tilde{F}_1) + \epsilon \frac{n^2}{m^2} B^1 \partial_{y^1}\tilde{F}_1 \right)\cdot \partial_{y^1} \tilde{F}_0
\end{eqnarray}
Next, examine the term containing $\partial_{y^0 y^0}\tilde{F}_0$. Using $\tilde{F}_0(y^0,y^1) = F_0(\frac{y^1 - a(y^0)}{\epsilon};R(y^0))$, the definition of $\tilde{F}_0$, to see that
\begin{eqnarray}\label{partial y^0 y^0 F_0 partial y^1 F_0 term}
 \int\limits_{-y_*^1}^{y_*^1} \partial_{y^0 y^0}\tilde{F}_0\cdot \partial_{y^1}\tilde{F}_0
 &=& -\frac{a''}{\epsilon^2} \int\limits_{-y_*^1}^{y_*^1} \partial_{y^1} F_0^2  + \frac{1}{\epsilon}\int\limits_{-y_*^1}^{y_*^1} \left[ (\frac{a'}{\epsilon})^2 \partial_{y^1 y^1}F_0 - 2 R' \frac{a'}{\epsilon} \partial_{y^1}\partial_R F_0 \right. \\
 & & \;\;\;\;\;\;\;\;\;\;\;\;\;\;\;\;\;\;\;\;\;\;\;\;\;\;\;\; \left. + R'' \partial_R F_0 + (R')^2 \partial_{RR} F_0 \right]\cdot \partial_{y^1} F_0 \nonumber
\end{eqnarray}
where again we've used the notation $\partial_{y^1}^\alpha \partial_R^\beta F_i = \partial_{y^1}^\alpha\partial_R^\beta F_i(\frac{y^1 - a}{\epsilon};R)$. We would like to obtain a bound for $\left| \frac{a''}{\epsilon}\right|$ in order to control $\underline{A}$. To do this, we will use (\ref{xi perp partial y^1 F_0} - \ref{partial y^0 y^0 F_0 partial y^1 F_0 term}) and isolate for the $\partial_{y^1} F_0^2$ term. We will then use this expression to obtain theorem \ref{bounded shift theorem}.

\begin{enumerate}
 \item Using (\ref{partial y^0 y^0 F_0 partial y^1 F_0 term}), Cauchy-Schwarz, and corollary \ref{partial F_i tilde estimate corollary} we have
 \[
   \frac{\left| a''\right|}{\epsilon} \int\limits_{\left| y^1\right| \leq y_*^1} \partial_{y^1} F_0^2 \lesssim \left|\; \int\limits_{\left| y^1\right| \leq y_*^1} \partial_{y^0 y^0} \tilde{F}_0\cdot \partial_{y^1}\tilde{F}_0 \right| + \underline{A}
 \]

 \item Using (\ref{S_0 partial y^1 F_0 term}), Cauchy-Schwarz, and corollary \ref{partial F_i tilde estimate corollary} we have
 \[
  \left|\; \int\limits_{\left| y^1\right| \leq y_*^1} \partial_{y^0 y^0} \tilde{F}_0\cdot \partial_{y^1} \tilde{F}_0\right| \lesssim \left| \;\int\limits_{\left| y^1\right| \leq y_*^1} \frac{n^2}{m^2} S_0\cdot \partial_{y^1} \tilde{F}_0\right| + (1 + \left| a''\right|)\underline{A}
 \]
 
 \item Using (\ref{partial y^0 y^0 xi partial y^1 F_0}), Cauchy-Schwarz, and corollary \ref{partial F_i tilde estimate corollary} we have
 \begin{align*}
   \left| \;\int\limits_{\left| y^1\right| \leq y_*^1} \frac{n^2}{m^2} S_0\cdot \partial_{y^1}\tilde{F}_0\right| 
   & \lesssim \left| \;\int\limits_{\left| y^1\right| \leq y_*^1} \partial_{y^0 y^0}\xi \cdot\partial_{y^1}\tilde{F}_0\right| + \left| \;\int\limits_{\left| y^1\right| \leq y_*^1} \frac{n^2}{m^2} B^\alpha \partial_\alpha \xi \cdot \partial_{y^1}\tilde{F}_0\right| \\
   & \; + \left| \; \int\limits_{\left| y^1\right| \leq y_*^1} \frac{n^2}{m^2} L_\epsilon(\tilde{F}_0,R) \xi \cdot \partial_{y^1} \tilde{F}_0\right| + \left| \; \int\limits_{\left| y^1\right| \leq y_*^1} \frac{n^2}{m^2} S_{-1}\cdot \partial_{y^1} \tilde{F_0}\right| + \left| \; \int\limits_{\left| y^1\right| \leq y_*^1} \frac{n^2}{m^2} N\cdot \partial_{y^1}\tilde{F}_0\right|
 \end{align*}
 We will estimate the $B^\alpha \partial_\alpha \xi$, $L_\epsilon(\tilde{F}_0,R)\xi$, $S_{-1}$, and $N$ terms separately.
 \begin{enumerate}
  \item \underline{$B^\alpha \partial_\alpha\xi$ term}: Recall the definitions for $B^0$ (\ref{B^0}), $B^1$ (\ref{B^1}), and $E$ (\ref{energy}). Using Cauchy-Schwarz, corollary \ref{partial F_i tilde estimate corollary}, and the boundedness of $\frac{n^2}{m^2} B^\alpha$ on $(0,y_*^0)\times (-y_*^1,y_*^1)$ we have
  \begin{gather*}
   \left| \; \int\limits_{\left| y^1\right|\leq y_*^1} \frac{n^2}{m^2} B^\alpha \partial_\alpha\xi\cdot\partial_{y^1}\tilde{F}_0\right| 
   \lesssim E^{1/2} \left\|\partial_{y^1} \tilde{F}_0\right\|_{L^2(\mathbb{R})} 
   \lesssim \frac{1}{\sqrt{\epsilon}} E^{1/2} \underline{A}
  \end{gather*}
  
  \item \underline{$L_\epsilon(\tilde{F}_0,R)\xi$ term}: Integrating by parts with respect to $y^1$ twice, using that $\operatorname{Hess}_\Phi W(\tilde{F}_0,R)$ is symmetric, using that $\partial_{y^1} \tilde{F}_0\in \ker(L_\epsilon(\tilde{F}_0,R)$, and using that on $(0,y_*^0)\times (-y_*^1,y_*^1)$, $\frac{n^2}{m^2}$ is bounded we have
  \begin{align*}
   \left| \; \int\limits_{\left| y^1\right|\leq y_*^1} \frac{n^2}{m^2} L_\epsilon(\tilde{F}_0,R)\xi\cdot\partial_{y^1}\tilde{F}_0\right| 
   \leq& \left|\; \int\limits_{\left| y^1\right|\leq y_*^1} \partial_{y^1 y^1} \left( \frac{n^2}{m^2}\right) \xi\cdot\partial_{y^1} \tilde{F}_0\right| + \left| \; \int\limits_{\left| y^1\right|\leq y_*^1} 2\partial_{y^1}\left(\frac{n^2}{m^2}\right) \xi\cdot \partial_{y^1 y^1}\tilde{F}_0 \right| \\
   & + \left| \left[\partial_{y^1}\xi\cdot\partial_{y^1}\tilde{F}_0 + \xi\cdot\partial_{y^1}\left(\frac{n^2}{m^2}\right) \partial_{y^1} \tilde{F}_0\right]_{-y_*^1}^{y_*^1}\right| \\
   \leq& \frac{1}{\sqrt{\epsilon}} E^{1/2} + \frac{1}{\epsilon^2} e^{-\alpha \frac{y_*^1}{\epsilon}} e^{\alpha \underline{A}}
  \end{align*}

  \item \underline{$S_{-1}$ term}: Using the definition of $S_{-1}$ (\ref{S_-1}), Cauchy-Schwarz, corollary \ref{partial F_i tilde estimate corollary}, using the fact that $B^1 = -H(R) + O(y^1)$, and the boundedness of $\frac{n^2}{m^2}$ on $(0,y_*^0)\times (-y_*^1,y_*^1)$ we have
  \begin{eqnarray*}
   \left| \; \int\limits_{\left| y^1\right|\leq y_*^1} S_{-1}\cdot \partial_{y^1} \tilde{F}_0\right|
   &=& \left| \; \int\limits_{\left| y^1\right|\leq y_*^1} \frac{n^2}{m^2} \partial_{y^1} \tilde{F}_0\cdot\left[B^1 + H(R)\right]\partial_{y^1} \tilde{F}_0 + \frac{n^2}{m} \frac{a}{\epsilon} \partial_R w(\tilde{F}_0,R) \cdot\partial_{y^1} \tilde{F}_0\right| \\
   &\lesssim& \int\limits_{-y_*^1}^{y_*^1} \left| y^1\right| \left| \partial_{y^1} \tilde{F}_0\right|^2 + \frac{\left| a\right|}{\epsilon} \left\| \partial_R w(\tilde{F}_0,R)\right\|_{L^2(\mathbb{R})} \left\| \partial_{y^1} \tilde{F}_0\right\|_{L^2(\mathbb{R})} \\
   &\lesssim& \underline{A}
  \end{eqnarray*}

 \item \underline{$N$ term}: To estimate this term we proceed as we did in the energy estimate when we estimated the $N\cdot \partial_{y^0} \xi$ term in lemma \ref{N estimate lemma}, where $N$ was defined in (\ref{N}). To obtain this estimate, we again use the identity $g(t) = g(0) + g'(0)t + \int\limits_0^1 (1-t) g''(t)$ to rewrite $N$ as 
 \begin{align*}
  N 
  &= -\frac{1}{\epsilon^2}\int\limits_0^1 (1-t) \frac{d^2}{dt^2} w(\tilde{F}_0 + \tilde{F}_\xi, R + t y^1 m)dt
 \end{align*}
 Thus,
 \begin{align*}
  \epsilon^2 \left| \int\limits_{\left| y^1\right| \leq y_*^1} \frac{n^2}{m^2} N\cdot \partial_{y^1}\tilde{F}_0\right|
  =& \; \epsilon^2 \left| \left< \frac{n^2}{m^2} N,\partial_{y^1}\tilde{F}_0\right>_{L^2(-y_*^1,y_*^1)}\right| \\
  \leq& \max\limits_{t\in [0,1]}\left| \left< \frac{n^2}{m^2} \frac{d^2}{dt^2} w(\tilde{F}_0 + \tilde{F}_\xi, R + t y^1 m) ,\partial_{y^1}\tilde{F}_0\right>_{L^2(-y_*^1,y_*^1)}\right| \\
  \leq& \max\limits_{t\in [0,1]} \left|\left< \left( (\tilde{F}_\xi)_\phi \partial_\phi + (\tilde{F}_\xi)_\sigma \partial_\sigma\right) \operatorname{Hess}_\Phi W(\tilde{F}_0 + t\tilde{F}_\xi,R + t y^1 m)\tilde{F}_\xi ,\partial_{y^1}\tilde{F}_0\right>_{L^2(-y_*^1,y_*^1)}\right| \\
  & + \; \max\limits_{t\in [0,1]} \left|\left< 2 y^1 m(y^0) \partial_R \operatorname{Hess}_\Phi W(\tilde{F}_0 + t\tilde{F}_\xi,R + t y^1 m) \tilde{F}_\xi ,\partial_{y^1}\tilde{F}_0\right>_{L^2(-y_*^1,y_*^1)}\right| \\
  & + \; \max\limits_{t\in [0,1]} \left|\left< (y^1 m)^2 \partial_{RR} w  ,\partial_{y^1}\tilde{F}_0\right>_{L^2(-y_*^1,y_*^1)}\right|
 \end{align*}
 and hence
 \begin{align*}
  & \left| \int\limits_{\left| y^1\right| \leq y_*^1} \frac{n^2}{m^2} N\cdot \partial_{y^1}\tilde{F}_0\right| \\
   \lesssim& \frac{1}{\epsilon^2} \left\| \tilde{F}_\xi\right\|_\infty \left\|\tilde{F}_\xi\right\|_{L^2(-y_*^1,y_*^1)} \left\| \partial_{y^1} \tilde{F}_0\right\|_{L^2(\mathbb{R})} + \frac{1}{\epsilon^2} \left\| \tilde{F}_\xi\right\|_{L^2(-y_*^1,y_*^1)} \left\| y^1 \partial_{y^1} \tilde{F}_0\right\|_{L^2(\mathbb{R})} + \frac{1}{\epsilon^2} \int\limits_{-y_*^1}^{y_*^1} (y^1)^2 \tilde{s}_0 \partial_{y^1} \tilde{s}_0 \\
  \lesssim& \frac{1}{\epsilon^{5/2}} (\epsilon + \left\|\xi\right\|_\infty)(\epsilon^{3/2} + \epsilon E^{1/2} + e^{-\alpha \frac{y_*^1 - a}{\epsilon}}) + \frac{1}{\epsilon^{3/2}}(\epsilon^{3/2} + \epsilon E^{1/2} + e^{-\alpha \frac{y_*^1 - a}{\epsilon}}) + 1
 \end{align*}
 We are left to estimate $\left\| \xi\right\|_\infty$. Using Gagliardo-Nirenberg we get that 
 \begin{eqnarray*}
  \left\| \xi\right\|_\infty
  &\lesssim& \left\|\xi\right\|_{L^2(-y_*^1, y_*^1)}^{1/2} \left\| \partial_{y^1} \xi\right\|_{L^2(-y_*^1, y_*^1)}^{1/2} \\
  &\lesssim& (\epsilon E^{1/2} + e^{-\alpha \frac{y_*^1 - a}{\epsilon}})^{1/2} E^{1/4}
 \end{eqnarray*}
 Thus, we have that 
 \begin{eqnarray*}
  \left| \; \int\limits_{-y_*^1}^{y_*^1} \frac{n^2}{m^2} N\cdot\partial_{y^1} \tilde{F}_0\right|
  &\lesssim& 1 + \frac{1}{\epsilon^{5/2}} (\epsilon + (\epsilon E^{1/2} + e^{-\alpha \frac{y_*^1 - a}{\epsilon}})^{1/2} E^{1/4})(\epsilon^{3/2} + \epsilon E^{1/2} + e^{-\alpha \frac{y_*^1 - a}{\epsilon}})
 \end{eqnarray*}
 \end{enumerate}
 
 \item Finally, recall the definition of $E$ (\ref{energy}). Using (\ref{partial y^0 y^0 xi partial y^1 F_0}), Cauchy-Schwarz, and corollary \ref{partial F_i tilde estimate corollary} we have
 \[
  \left| \; \int\limits_{\left| y^1\right| \leq y_*^1} \partial_{y^0 y^0}\xi\cdot \partial_{y^1}\tilde{F}_1\right| \lesssim \frac{1}{\sqrt{\epsilon}}(E + \frac{1}{\epsilon} e^{-\alpha \frac{y_*^1 - a}{\epsilon}}) \left[1 + \frac{\left| a''\right|}{\epsilon}\right]\underline{A}
 \] 
\end{enumerate}

Combining the estimates obtained in steps 1-4 we get that
\begin{eqnarray*}
 \frac{\left|a''\right|}{\epsilon}
 &\lesssim& \frac{\left| a'' \right|}{\epsilon} \left[ \epsilon \underline{A} + \sqrt{\epsilon} \underline{A} \sqrt{E} + e^{-\alpha \frac{y_*^1 - a}{\epsilon}} \right] + \left[ \frac{1}{\sqrt{\epsilon}^{5}} (\epsilon\underline{A} + (\epsilon \sqrt{E} + e^{-\alpha \frac{y_*^1 - a}{\epsilon}})^{1/2} \sqrt[4]{E})(\sqrt{\epsilon}^{3} + \epsilon \sqrt{E} + e^{-\alpha \frac{y_*^1 - a}{\epsilon}}) \right]
\end{eqnarray*}
as desired. This concludes the proof of theorem \ref{bounded shift theorem}.
\begin{flushright}
 $\Box$
\end{flushright}

\begin{appendices}
\section{Formal Asymptotics}\label{formal asymptotics}
Let $\eta$ be the Minkowski metric on $\mathbb{R}^{1+n}$ and let $\Gamma\subset (\mathbb{R}^{1+n},\eta)$ be an $n$-dimensional time-like surface in space-time. Suppose that $\Gamma$ is parameterized by some map $H:\Omega\subset\mathbb{R}^n\rightarrow \mathbb{R}^{1+n}$. Define a new coordinate system $(y^\tau,y^\nu)\in\mathbb{R}^n\times \mathbb{R}$, called Minkowski normal coordinates, as $$(t,x) = \psi(y^\tau,y^\nu) = H(y^\tau) + y^\nu \nu(y^\tau)$$ where $\nu(y^\tau) \perp_\eta \partial_{y^\tau} H(y^\tau)$ and $\left| \nu(y^\tau) \right|_\eta = 1$. We call $y^\tau\in\mathbb{R}^{n}$ ``tangential coordinates'' and $y^\nu\in\mathbb{R}$ the ``normal coordinate''. Note that this coordinate system may only be well defined on a neighbourhood $\mathcal{N}$ of $\Gamma$.
\bigskip

Recall that we want to find solutions of (\ref{equations of motion of interface with a current}) so that $\phi$ has an interface and so that $\sigma$ is exponentially small except near the interface of $\phi$. Based on \cite{jerrard2011semilinear}, we expect that for suitable $\Gamma$, $\theta:\Omega\rightarrow\mathbb{R}$, and $\Phi_0 := (\phi_0,\sigma_0):\mathbb{R}\rightarrow\mathbb{R}^2$ there exists a solution with these characteristics of the form
\begin{equation}\label{approx solution 1}
 \Phi(y^\tau,y^\nu) \approx \left(\begin{array}{c} \phi_0(\frac{y^\nu}{\epsilon}) \\ e^{\frac{i}{\epsilon} \theta(y^\tau)} \sigma_0(\frac{y^\nu}{\epsilon}) \end{array}\right)
\end{equation}
where $\gamma_{ij} := \eta_{\alpha\beta} \partial_i H^\alpha \partial_j H^\beta$ is the induced metric on the surface $\Gamma$ (latin indices range over the tangential coordinates and Greek indices will range over both tangential and normal coordinates).
\bigskip

We will now carry out a formal asymptotic analysis to find $\Phi_0$ so that $\phi_0$ has an interface and to find $\Gamma$ and $\theta$ for which we expect (\ref{approx solution 1}) to hold. To do this, we will expand the action integral associated to (\ref{equations of motion of interface with a current}) about the right hand side of (\ref{approx solution 1}). From this expansion, we obtain an \textbf{effective action}. We will then make a choice for the profile $\Phi_0$ and for this choice of $\Phi_0$, we expect, heuristically, that the correction terms coming from expanding the action about the right hand side of (\ref{approx solution 1}) will be of lower order when $\Gamma$ and $\theta$ are critical points of the effective action.
\bigskip

The Lagrangian associated to (\ref{equations of motion of interface with a current}) in Minkowski normal coordinates is
\begin{equation}\label{Lagrangian}
 \mathcal{L} := \frac{1}{2} g^{\alpha\beta} \partial_\alpha \phi \partial_\beta \phi + \frac{1}{2} g^{\alpha\beta} \overline{\partial_\alpha \sigma} \partial_\beta \sigma + \frac{1}{\epsilon^2} V(\phi,\sigma)
\end{equation}
where $g_{\alpha\beta}:= \eta_{\lambda\omega} \partial_\alpha \psi^\lambda \partial_\beta \psi^\omega$ is the Minkowski metric in normal coordinates. Note that
\[
g_{\alpha\beta} = \left(\begin{array}{c:c} \gamma_{ij}  & 0 \\ \hdashline 0 & 1  \end{array}\right) +  (y^\nu)^2 \left(\begin{array}{c:c} \eta_{\lambda\omega} \partial_i \nu^\lambda \partial_j \nu^\omega & 0 \\ \hdashline 0 & 0  \end{array}\right)
\]
For $\xi = (\xi_\phi,\xi_\sigma):\mathbb{R}^{1+n}\rightarrow\mathbb{R}\times\mathbb{C}$, we plug
\begin{gather*}
 \phi = \phi_0(\frac{y^\nu}{\epsilon}) + \xi \;\;\;\;\; \text{ and } \;\;\;\;\; \sigma = e^{\frac{i}{\epsilon}\theta(y^\tau)} \left[\sigma_0(\frac{y^\nu}{\epsilon}) + \xi_\sigma\right]
\end{gather*}
into the action integral to get
\begin{equation}\label{action integral ansatz}
 S(\Phi) = \frac{1}{\epsilon^2} \int \left\{ \frac{1}{2} (\Phi_0'(\frac{y^\nu}{\epsilon}))^2 + V(\Phi_0)(\frac{y^\nu}{\epsilon}) + \frac{1}{2} \gamma^{ij}\partial_i\theta\partial_j\theta \sigma_0(\frac{y^\nu}{\epsilon})^2 \right\} \sqrt{-\gamma(y^\tau)} dy^\tau dy^\nu + \text{other terms}
\end{equation}
The effective action we obtain from this expansion is
\begin{equation}\label{effective action}
 \tilde{S} := \int \left\{ \frac{1}{2} (\Phi_0'(\frac{y^\nu}{\epsilon}))^2 + V(\Phi_0)(\frac{y^\nu}{\epsilon}) + \frac{1}{2} \gamma^{ij}\partial_i\theta\partial_j\theta \sigma_0(\frac{y^\nu}{\epsilon})^2 \right\} \sqrt{-\gamma(y^\tau)} dy^\tau dy^\nu
\end{equation}
\bigskip

Consider the $\frac{1}{\epsilon^2}$ term. It is natural to choose $\Phi_0$ so that in transverse directions to $\Gamma$, $\Phi_0$ is energy minimizing and so that $\phi_0$ has an interface. To this end, suppose for $\rho\in\mathbb{R}$, $F = (f,s)(\cdot;\rho)$ satisfies the minimization problem
\begin{gather}
 \mu(\rho) := \inf\limits_{(f,s)\in\mathcal{A}} \int\left\{ \frac{1}{2} \right|(f',s')\left|^2 + V(f,s) + \frac{1}{2} \rho s^2\right\} dy^\nu \label{transverse energy min} \\
 \mathcal{A} := \left\{ (f,s) \in C^1 \; : \; \lim\limits_{y^\nu\rightarrow \pm \infty} f(y^\nu) = \pm 1  \right\} \label{boundary conditions}
\end{gather}
In this case, the boundary conditions imposed results in $f$ having an interface. Furthermore, for suitable potentials $V$, $s$ is exponentially small except near the interface of $f$. We pick $\Phi_0 = (f,s)(\cdot;\zeta)$, where $\zeta(y^\tau) := \gamma^{ij} \partial_i\theta\partial_j\theta$. \textit{Important}: The natural choice of profile $\Phi_0$ actually depends on $\zeta$. That is, in contrast to our initial hypothesis (\ref{approx solution 1}), we expect that there should exist a solution to (\ref{equations of motio7n of interface with a current}) satisfying
\begin{equation}\label{approx solution 2}
 \Phi \approx \left(\begin{array}{c} \phi_0(\frac{y^\nu}{\epsilon};\zeta(y^\tau)) \\ e^{\frac{i}{\epsilon}\theta(y^\tau)} \sigma_0(\frac{y^\nu}{\epsilon};\zeta(y^\tau)) \end{array}\right)
\end{equation}
for $\Phi_0 = \Phi_0(\cdot;\zeta)$ minimizing (\ref{transverse energy min}) and for suitable $\Gamma$ and $\theta$.
\bigskip

For this choice of $\Phi_0$, the effective action becomes
\begin{equation}\label{effective action 2}
 \tilde{S}(H,\theta) = \int \mu(\zeta) \sqrt{-\gamma} dy^\tau
\end{equation}
Heuristically, we expect that when $\theta$ and $H$ are critical points of $\tilde{S}$, then $\xi$ will be of lower order than the right hand side of (\ref{approx solution 1}). That is, for $\theta$ and $H$ satisfying the nonlinear, coupled system
\begin{align}\label{EL for theta}
 0 = \frac{\delta \tilde{S}}{\delta \theta} = -2 \partial_j \left( \mu'(\zeta) \sqrt{-\gamma} \gamma^{ij} \partial_i \theta \right)
\end{align}
\begin{align}\label{EL for H}
 0 = \frac{\delta S}{\delta H} = -\eta_{\alpha\beta} \partial_j \left( \mu(\zeta) \sqrt{-\gamma} \gamma^{ij} \partial_i H^\alpha \right) + 2 \eta_{\alpha\beta} \partial_j \left( \mu'(\zeta) \sqrt{-\gamma} \gamma^{ik} \gamma^{lj} \partial_k \theta \partial_l \theta \partial_i H^\beta\right)
\end{align}
then $\Phi_0(\frac{y^\nu}{\epsilon};\zeta)$ should be a good approximate solution. The coupled system for $\theta$ and $H$ should be a hyperbolic system, but this isn't completely clear. By expanding (\ref{EL for H}) and taking its inner product with $\nu^\beta$, we can rewrite this system as
\begin{gather}
 \Box_\Gamma \theta = -\gamma(\nabla_\tau\log\left[\mu'(\zeta)\right],\nabla_\tau\theta) \label{geometric EL for theta} \\
 \text{mean curvature of } \Gamma = 2 \frac{\mu'(\zeta)}{\mu(\zeta)} \mathrm{I\!I}\left(\nabla_\tau \theta,\nabla_\tau\theta\right) \label{geometric EL for H}
\end{gather}
where $\mathrm{I\!I}$ is the second fundamental form of $\Gamma$. From (\ref{geometric EL for H}) we find a nice geometric relation between the surface about which our approximate solution is concentrated and the phase of $\sigma_0$.

\end{appendices}

\bibliographystyle{plain}
\bibliography{Bibliography.bib}

\begin{thebibliography}{10}

\bibitem{aftalion2012vortex}
A.~Aftalion, P.~Mason, and J.~Wei.
\newblock Vortex-peak interaction and lattice shape in rotating two-component
  bose-einstein condensates.
\newblock {\em Physical Review A}, 85(3):033614, 2012.

\bibitem{aftalion2015thomas-fermi}
A.~Aftalion, B.~Noris, and C.~Sourdis.
\newblock Thomas-fermi approximation for coexisting two component bose-einstein
  condensates and nonexistence of vortices for small rotation.
\newblock {\em Comm. Math. Phys.}, 336(2):509--579, 2015.

\bibitem{alama2006fractional}
S.~Alama and L.~Bronsard.
\newblock Fractional degree vortices for a spinor ginzburg-landau model.
\newblock {\em Communications in Contemporary Mathematics}, 8(03):355--380,
  2006.

\bibitem{alama2015domain}
S.~Alama, L.~Bronsard, A.~Contreras, and D.~Pelinovsky.
\newblock Domain walls in the coupled {G}ross-{P}itaevskii equations.
\newblock {\em Arch. Ration. Mech. Anal.}, 215(2):579--610, 2015.

\bibitem{alama2009structure}
S.~Alama, L.~Bronsard, and P.~Mironescu.
\newblock On the structure of fractional degree vortices in a spinor
  ginzburg-landau model.
\newblock {\em Journal of Functional Analysis}, 256(4):1118--1136, 2009.

\bibitem{alama2012compound}
S.~Alama, L.~Bronsard, and P.~Mironescu.
\newblock On compound vortices in a two-component {G}inzburg-{L}andau
  functional.
\newblock {\em Indiana Univ. Math. J.}, 61(5):1861--1909, 2012.

\bibitem{alikakos2008connection}
N.~Alikakos and G.~Fusco.
\newblock On the connection problem for potentials with several global minima.
\newblock {\em Indiana Univ. Math. J}, 57(4):1871--1906, 2008.

\bibitem{coddington1955theorem}
E.~Coddington and N.~Levinson.
\newblock {\em Theorem of Ordinary Differential Equations}.
\newblock McGraw-Hill, 1955.

\bibitem{czubak2015topological}
M.~Czubak and R.~Jerrard.
\newblock Topological defects in the abelian {H}iggs model.
\newblock {\em Discrete Contin. Dyn. Syst.}, 35(5):1933--1968, 2015.

\bibitem{galvao2015accelerating}
B.~Galv{\~a}o-Sousa and R.~Jerrard.
\newblock Accelerating fronts in semilinear wave equations.
\newblock {\em Rend. Circ. Mat. Palermo (2)}, 64:117--148, 2015.

\bibitem{gustafson2006effective}
S.~Gustafson and I.~M. Sigal.
\newblock Effective dynamics of magnetic vortices.
\newblock {\em Adv. Math.}, 199(2):448--498, 2006.

\bibitem{isoshima2002axisymmetric}
T.~Isoshima and K.~Machida.
\newblock Axisymmetric vortices in spinor bose-einstein condensates under
  rotation.
\newblock {\em Physical Review A}, 66(2):023602, 2002.

\bibitem{jerrard1999vortex}
R.~Jerrard.
\newblock Vortex dynamics for the {G}inzburg-{L}andau wave equation.
\newblock {\em Calc. Var. Partial Differential Equations}, 9(1):1--30, 1999.

\bibitem{jerrard2011semilinear}
R.~Jerrard.
\newblock Defects in semilinear wave equations and timelike minimal surfaces in
  minkowski space.
\newblock {\em Anal. PDE}, 4(2):285--340, 2011.

\bibitem{jerrard2015dynamics}
R.~Jerrard.
\newblock Dynamics of topological defects in nonlinear field theories.
\newblock {\em Adv. Stud. Pure. Math}, 67:157--224, 2015.

\bibitem{knigavko1998spontaneous}
A.~Knigavko and B.~Rosenstein.
\newblock Spontaneous vortex state and ferromagnetic behavior of type-ii p-wave
  superconductors.
\newblock {\em Physical Review B}, 58(14):9354, 1998.

\bibitem{lin1999vortex}
F.~Lin.
\newblock Vortex dynamics for the nonlinear wave equation.
\newblock {\em Comm. Pure Appl. Math.}, 52(6):737--761, 1999.

\bibitem{mason2011classification}
P.~Mason and A.~Aftalion.
\newblock Classification of the ground states and topological defects in a
  rotating two-component bose-einstein condensate.
\newblock {\em Physical Review A}, 84(3):033611, 2011.

\bibitem{nielsen1973vortex}
H.~Nielsen and P.~Olesen.
\newblock Vortex-line models for dual strings.
\newblock {\em Nuclear Physics B}, 61:45--61, 1973.

\bibitem{reed1980functional}
M.~Reed and B.~Simon.
\newblock {\em Methods of Modern Mathematical Physics}, volume~1.
\newblock Academic Press, 1980.

\bibitem{shatah1998geometric}
J.~Shatah and M.~Struwe.
\newblock {\em Geometric Wave Equations}, volume~2.
\newblock American Mathematical Soc., 1998.

\bibitem{stuart2004geodesic}
D.~Stuart.
\newblock The geodesic hypothesis and non-topological solitons on
  pseudo-{R}iemannian manifolds.
\newblock {\em Ann. Sci. \'Ecole Norm. Sup. (4)}, 37(2):312--362, 2004.

\bibitem{stuart2004geodesics}
D.~Stuart.
\newblock Geodesics and the einstein nonlinear wave system.
\newblock {\em Journal de math{\'e}matiques pures et appliqu{\'e}es},
  83(5):541--587, 2004.

\bibitem{tao2006nonlinear}
T.~Tao.
\newblock {\em Nonlinear Dispersive Equations: Local and Global Analysis},
  volume 106.
\newblock American Mathematical Soc., 2006.

\bibitem{vilenkin2000cosmic}
A.~Vilenkin and P.~Shellard.
\newblock {\em Cosmic Strings and Other Topological Defects}.
\newblock Cambridge University Press, 2000.

\bibitem{witten1985superconducting}
E.~Witten.
\newblock Superconducting strings.
\newblock {\em Nuclear Physics B}, 249(4):557--592, 1985.

\end{thebibliography}

\end{document}